\documentclass[preprint,11pt]{elsarticle}

\usepackage{siunitx}
\usepackage{booktabs}
\usepackage{amsthm}
\usepackage{subcaption}

\usepackage{amsmath,amssymb,bm}
\usepackage{accents}

\biboptions{sort}
\usepackage[usenames,dvipsnames,svgnames,table]{xcolor}

\usepackage{color}

\usepackage{tikz}
\usetikzlibrary{plotmarks}
\usetikzlibrary{positioning}
\usetikzlibrary{decorations.pathreplacing}
\usetikzlibrary{math}

\usepackage{graphicx}
\usepackage{calc,fancyvrb}
\usepackage[margin=1.in]{geometry}

\usepackage{algorithm,algpseudocode}
\usepackage{algorithmicx}
\algrenewcommand\alglinenumber[1]{\footnotesize #1:} 
\newcommand{\algFontSize}{\small}

\usepackage{pgfplots}
\usepackage{graphicx}
\pgfplotsset{compat=1.16}
\usepackage[bookmarks=true,colorlinks=true,linkcolor=blue]{hyperref}

\newtheorem{theorem}{Theorem}
\newtheorem{proposition}{Proposition}[section]
\newtheorem{lemma}{Lemma}[section]
\newtheorem{definition}{Definition}[section]

\newtheorem{remark}{Remark}[section]

\numberwithin{equation}{section}

\begin{document}

\begin{frontmatter}
 \title{Local-basis Difference Potentials Method for elliptic PDEs in complex geometry}

\author{Qing~Xia}
\ead{qingx@kth.se}

\address{Department of Mathematics, KTH Royal Institute of Technology, Stockholm, Sweden}

\begin{abstract}
We develop efficient and high-order accurate finite difference methods for elliptic partial differential equations in complex geometry in the Difference Potentials framework. The main novelty of the developed schemes is the use of local basis functions defined at near-boundary grid points. The use of local basis functions allow unified numerical treatment of (i) explicitly and implicitly defined geometry; (ii) geometry of more complicated shapes, such as those with corners, multi-connected domain, etc; and (iii) different types of boundary conditions. This geometrically flexible approach is complementary to the classical difference potentials method using global basis functions, especially in the case where a large number of global basis functions are needed to resolve the boundary, or where the optimal global basis functions are difficult to obtain. Fast Poisson solvers based on FFT are employed for standard centered finite difference stencils regardless of the designed order of accuracy. Proofs of convergence of difference potentials in maximum norm are outlined both theoretically and numerically. 
\end{abstract}


\end{frontmatter}



\section{Introduction}

In this work, we develop efficient and high-order accurate numerical schemes for elliptic partial differential equations in complex geometry. The designed methods are based on local basis functions and finite difference methods, which is orchestrated in the difference potentials framework. Details of difference potentials can be found for example, in the monograph \cite{ryaben2001method}. Difference Potentials Method (DPM) is also based on the formulation of boundary integrals, but differs with boundary integral methods in using finite difference / finite volume / finite element discretizations to approximate singular integrals that arise in the solutions of elliptic partial differential equations, which avoids tackling singular integrals with special quadratures directly. DPM is similar to boundary integral methods in dimension reduction, in the sense that it reduces the spatial discretizations to algebraic equations defined only at grid points near the boundary, but requires no knowledge of the fundamental solutions. DPM has been studied and proposed for various types of model problems, including but not limited to, interface problems of elliptic or parabolic types \cite{epshteyn2015high,Albright_2017_parabolic,Albright_2017_elliptic,Ludvigsson_2018}, Helmholtz equations \cite{britt2015computation,medvinsky2016solving,north2022non}, wave problems in 2D and 3D \cite{Britt_2018,Petropavlovsky_2018,Medvinsky_2013,medvinsky2021solution,petropavlovsky20223d}, chemotaxis models in 2D and 3D \cite{Epshteyn_2012,MR3954435}, bulk-surface problems in 3D \cite{epshteyn2020difference}, etc. 

In the classical difference potentials framework, an important assumption is that the boundary or interface data, i.e. $\partial^ku/\partial n^k$ for $k=0,1,2,\dots$ with $n$ denoting outward unit normal vector, can be expressed by finite sums of global basis functions defined on the boundary or interface. Such spectral representations introduce excellent efficiency in the Difference Potentials Method, if a relatively small number of global basis functions can resolve the boundary or interface data sufficiently. The small number of global basis functions leads to a low-rank boundary equations formulated at grid points near the boundary, which means they can be solved efficiently. This condition on the number of global basis functions normally is satisfied, when the underlying geometry and the boundary/interface data are sufficiently smooth. In practice, good choices of global basis functions can reduce the total number of basis functions needed significantly. For instance, Fourier basis functions are the optimal candidates for circles/spheres; piecewise-smooth Fourier basis functions or Chebyshev polynomials work best for piecewise smooth geometry \cite{Magura_2017}. It is still an open and challenging question that how the optimal basis functions can be constructed accurately and efficiently for geometry of arbitrary shapes. In theory, one possible approach is to approximate eigenfunctions of the Laplace-Beltrami operator on the corresponding geometry, but is difficult to put this approach into practice. In addition, this ``small-number'' condition will not be satisfied even for spheres, if functions on the sphere need to be resolved by a large number of spherical harmonics especially when singularities occur on the surface, even though we already have the optimal basis functions on the sphere, namely spherical harmonics.


The local-basis-function approach we develop in this manuscript aims at complementing the classical approach of using global basis functions in the difference potentials framework, especially when the strength of global-basis approach falls short. 
In comparison, the local-basis-function approach also inherits the strengths of the classical global-basis-function approach in the difference potentials framework. For example, the introduction of computationally-simple auxiliary domains allow the employment of fast Poisson solvers, even though the underlying geometry is irregular.  The dimension reduction property is preserved -- interior discrete PDEs are reduced to boundary equations defined at grid points near the boundary similarly in both the global approach. The exact definitions of the basis functions are different: the global basis functions are defined on the entire boundary or interface, while the local basis functions are defined at grid points near the boundary or interface, which have compact supports like in the finite element approach. Thus each local basis function only interacts with a small part of curve of surface and requires no strict global regularity. For multi-connected domains, the local-basis-function approach only needs information of interior and exterior points, while the global approach needs to track each boundary or interface.
In terms of boundary conditions (BCs), the global-basis-function approach has flexibility in handling different kinds of BCs, thanks to the assumption on the spectral representation of boundary data. The local-basis-function approach is also flexible in handling BCs using a collocation approach, but is with drawback of losing one order of accuracy for Neumann boundary conditions, which we will continue to address in future work. Another feature of the local-basis-function approach is the unified treatment of explicit and implicit representations of the geometry. Implicit geometry was first treated using the difference potentials framework in \cite{Ludvigsson_2018}, where non-parameterized curves were converted into parametric ones and the rest follows the algorithms for explicit geometry. The local-basis approach developed here only requires knowledge of interior and exterior domains, which is supplied by level set functions. Thus the local-basis-function approach is insensitive to the exact representation of curves/surfaces, while the global-basis-function approach requires parametric forms.

We should also note that the numerical treatment of elliptic partial differential equations in complex geometry has an enormously extensive literature both in single and composite domains. In general, the methods in the current literature can be categorized into two main groups, those based on boundary integral formulations and those based on spatial discretizations with corrections at boundary or interface. In the first category, we have Mayo's method \cite{mayo1984fast}, boundary integral methods (BIM) of different flavors \cite{bystricky2021accurate,helsing2009integral,kublik2013implicit,zhong2018implicit,ying2007kernel,ying2013kernel,xie2020fourth,cao2022kernel}, boundary element methods (BEM) \cite{steinbach2007numerical}. In the second group, an incomplete list gives immersed boundary method \cite{peskin2002immersed}, immersed interface methods \cite{leveque1997immersed}, immersed finite element method \cite{zhang2004immersed,guo2020immersed}, cut finite element method \cite{hansbo2002unfitted,jonsson2017cut,burman2022posteriori,larson2020stabilization}, etc. The difference potentials framework (both local and global approaches) shares the simplicity with immersed mesh approaches and the dimension reduction property with boundary integral formulations of BIM and BEM. The main novelty of our developed schemes in this work is using local basis functions defined at grid points only near the boundary to allow geometric flexibility in the cut-FEM fashion, yet with no need of stabilization for small cut cells. Fast Poisson solvers with finite difference schemes on uniform Cartesian meshes are also designed. We showed that fast Poisson solvers based on FFT can be used for wide centered finite difference schemes can be applied for the auxiliary problem, regardless of the desired order of accuracy. In the literature of difference potentials, only fourth order compact stencils were used for fast solvers, for example in \cite{Medvinsky_2013}. The techniques of local basis functions developed in this work can be seamlessly applied to time-dependent problems, if implicit time stepping schemes are employed in the discretizations.

The rest of the manuscript is organized as follows. Section \ref{sec:scheme} outlines details of the designed numerical schemes with emphasis on the constructions of difference potentials using local basis functions. Some justifications of the convergence of difference potentials are presented in Section \ref{sec:convergence}. Section \ref{sec:numerics} validates the developed schemes using different setups, and some concluding remarks and discussions on future directions are in Section \ref{sec:conclusion}.

\section{Numerical schemes}\label{sec:scheme}

We consider the following elliptic equations in an arbitrary domain $\Omega\subset\mathbb{R}^d$,
\begin{align}
    Lu = f, \quad \mbox{in} \;\Omega\label{eqn:pde}
\end{align}
subject to some boundary condition at the boundary $\Gamma:=\partial \Omega$
\begin{align}
    \ell u  = g, \quad \mbox{in} \; \Gamma\label{eqn:bc}
\end{align}
The operator $L$ can be any elliptic operator, but in this manuscript we will focus on the form $Lu:=\Delta u -\sigma u$, which may represent Poisson's equation, Helmholtz equations, or discretized equations for heat equations, reaction-diffusion, convection diffusion equations, wave equations, etc., when implicit time stepping is employed for temporal discretization. 

\subsection{One dimension}

In this section, we will illustrate the main steps of the local-basis difference potentials approach in one dimension, and consider the equation \eqref{eqn:pde} in the domain $\Omega=[a,b]$, with Dirichlet boundary condition
\begin{align}\label{eqn:bc_1d}
u(a) = u_a,\quad u(b) = u_b.
\end{align}

\begin{figure}[H]
\centering
\begin{tikzpicture}[scale=1.5]
\draw[thick,|-|] (-1.5,0) -- (-0.25,0);
\draw[thick] (-0.25,0) -- (4.25,0);
\draw[thick,|-|] (4.25,0) -- (5.5,0);
\filldraw (0.0,0) circle (1pt);
\filldraw (0.5,0) circle (1pt);
\filldraw (1.0,0) circle (1pt);
\filldraw (1.5,0) circle (1pt);
\filldraw (2.0,0) circle (1pt);
\filldraw (2.5,0) circle (1pt);
\filldraw (3.0,0) circle (1pt);
\filldraw (3.5,0) circle (1pt);
\filldraw (4.0,0) circle (1pt);

\draw (-1.0,0) circle (1pt);
\draw (-0.5,0) circle (1pt);

\draw (-0.5,0) +(-1.5pt,-1.5pt) rectangle +(1.5pt,1.5pt);
\draw (0,0) +(-1.5pt,-1.5pt) rectangle +(1.5pt,1.5pt);

\draw (4.5,0) circle (1pt);
\draw (5.0,0) circle (1pt);

\draw (4.5,0) +(-1.5pt,-1.5pt) rectangle +(1.5pt,1.5pt);
\draw (4,0) +(-1.5pt,-1.5pt) rectangle +(1.5pt,1.5pt);

\draw (-0.25,0.1) node[above] {$a$};
\draw (4.25,0.1) node[above] {$b$};
\draw (-1.5,0.1) node[above] {$a-\ell$};
\draw (5.5,0.1) node[above] {$b+\ell$};

\filldraw (1.0,-0.5) circle (1pt);
\draw (1.3,-0.5) node {$M^+$};
\draw (2,-0.5) circle (1pt);
\draw (2.3,-0.5) node {$M^-$};
\draw (3,-0.5) +(-1.5pt,-1.5pt) rectangle +(1.5pt,1.5pt);
\draw (3.3,-0.5) node {$\gamma$};
\end{tikzpicture}
\caption{An example of point sets with the second order accurate discretization}
\label{fig:1d-example-O2}
\end{figure}
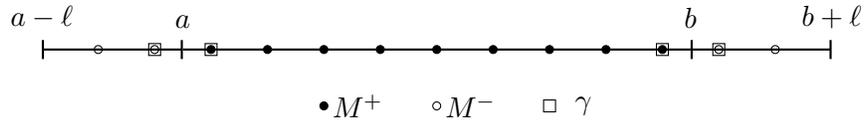

The first step in the difference potentials framework is to introduce an auxiliary domain $\Omega_0=[a-\ell,b+\ell]$ that embeds $\Omega$, then we discretize $\Omega_0$ using a uniform mesh with grid points $x_j$ ($j=0,1,\dots, N$). We will use centered finite difference stencils
\begin{align}\label{eqn:1d_stencil}
\mathcal{N}_{j}^{n} = \bigcup_{l=0}^n x_{j\pm l}
\end{align}
with $n=1$ for second order accuracy and $n=2$ for fourth order accuracy, etc.
The actual boundary points $x=a$ and $x=b$ may or may not align with a grid point. The choice of $\ell$ is also arbitrary, but should be comparably small for the length $b-a$. Next we define the point sets $M^\pm$, $M^0$, $N^{\pm},N^0$ and $\gamma^{\pm}$ needed in the schemes as follows.
\begin{definition}{(Point sets)}\label{def:point_sets} 
\hfill
\begin{itemize}
    \item $M^0= \left\{x_j\mid x_j\in\Omega^0\right\}$,\quad $M^+=M^0\cap\Omega$,\quad $M^-=M^0\backslash M^+$,

    \item $N^0=\left\{\bigcup_{j}\mathcal{N}_{j}^{n}\mid x_j\in M^0\right\}$,\quad $N^\pm=\left\{\bigcup_{j}\mathcal{N}_{j}^{n}\mid x_j\in M^\pm\right\}$,

    \item $\gamma=N^+\cap N^-$,\quad $\gamma^{\pm}=M^\pm\cap\gamma$.
\end{itemize}
\end{definition}
A graphical example of the defined point sets in the second order accurate case can be found in Figure~\ref{fig:1d-example-O2}. 

Assume the $n$-th order accurate finite difference discretization of \eqref{eqn:pde} is given by
\begin{align}
L_{nh}u_j = f_j,\quad x_j\in M^+. \label{eqn:discrete-eq}
\end{align}
To solve \eqref{eqn:discrete-eq} together with the boundary condition \eqref{eqn:bc}, the traditional finite difference approach is to either modify the stencil near the boundary or use interpolations between boundary data and the unknown interior nodes, which either lacks efficiency or accuracy. The difference potentials framework avoids such problem by using the auxiliary problems as defined below.
\begin{definition}{(Auxiliary problem)} \label{def:aux_prob}
We refer to solving the following discrete problem on $M^0$ as solving the auxiliary problem.
\begin{subequations}
\begin{align}
L_{nh}u_j &= q_j,\quad x_j \in M^0\\
\ell(u_j) &= p_j,\quad x_j\in N^0\backslash M^0
\end{align}
\end{subequations}
where $q_j$ is a grid function defined on $M^0$ and $p_j$ is some grid function on $N^0\backslash M^0$.
\end{definition}

\begin{remark}
The boundary condition of the auxiliary problem is arbitrary as long as it is well-posed. We will use this convenience to design fast solvers for the auxiliary problem on uniform meshes, which works for arbitrary choice of order of accuracy. 
\end{remark}

Next we define how to compute the particular solution and the difference potentials using the auxiliary problem defined in Definition~\ref{def:aux_prob} with some particular choices of $q_j$.
\begin{itemize}
    \item The particular solution $G_{nh}f$ is the solution to the auxiliary problem with the right hand side is
        \begin{align}\label{eqn:ps_rhs}
        q_j = \left\{
                \begin{array}{cc}
                f_j,& x_j\in M^+,\\
                0,& x_j\in M^-,
                \end{array}
                \right.
        \end{align}
        or in shorthand notation: $q=X_{M^+}f_j$, where $X_{M^+}$ denotes a characteristic function for set $M^+$. Then we denote the particular solution $G_{nh}f = G_{nh}X_{M^+}f_j$

    \item The difference potential $P_{N^+\gamma}v_\gamma$ associated with $v_\gamma$ is computed using the auxiliary problem with the following right hand side
        \begin{align}\label{eqn:dp_rhs}
        q_j = \left\{
                \begin{array}{cc}
                0,& x_j\in M^+,\\
                L_{nh}v_\gamma,& x_j\in M^-.
                \end{array}
                \right.
        \end{align}
        where $v_\gamma$ is the grid function defined on the point set $\gamma$ and we assume zero extension of $v_\gamma$ from $\gamma$ to $N^0$ when the operator $L_{nh}$ defined on $M^0$ is operates on $v_\gamma$. Similarly, we can denote $q=X_{M^-}L_{nh}v_\gamma$ and the difference potentials $P_{N^+\gamma}v_\gamma=G_{nh}X_{M^-}L_{nh}v_\gamma$.
\end{itemize}

The difference potential $P_{N^+\gamma}v_\gamma$ can be computed with any given density $v_\gamma$ using the auxiliary problem with the right hand side given in \eqref{eqn:dp_rhs}. There is only one density, which we denote as $u_\gamma$, that corresponds to the solution of \eqref{eqn:discrete-eq}, given the boundary condition \eqref{eqn:bc}, or in particular the boundary condition in 1D \eqref{eqn:bc_1d}. To this end, below we introduce a key theorem in the Difference Potentials framework, which prescribes how to obtain the desired density $u_\gamma$.

\begin{theorem}\label{thm:bep}
The grid function $u_\gamma$ is a trace to the solution of the discrete problem \eqref{eqn:discrete-eq}, if and only if the following boundary equations with projections hold:
\begin{align}\label{eqn:bep}
u_\gamma - P_\gamma u_\gamma = G_{nh}f_\gamma
\end{align}
where $P_\gamma u_\gamma:=Tr_\gamma P_{N^+\gamma}u_\gamma$ and $G_{nh}f_\gamma:=Tr_\gamma G_{nh}f$. The trace operator $Tr_\gamma$ denotes restriction to the point set $\gamma$.
\end{theorem}

Theorem~\ref{thm:bep} is a classical result in the difference potentials framework. The proof can be found in many references, including \cite{ryaben2001method}, \cite[Section 5.1]{MR3954435}, etc. An important implication of Theorem~\ref{thm:bep} is that the boundary equations \eqref{eqn:bep} is equivalent to the discrete PDE \eqref{eqn:discrete-eq} in $M^+$.

\begin{remark}
The difference potentials operator $P_{\gamma}$ is of dimension $|\gamma|\times|\gamma|$, and the entries in $P_{\gamma}$ can be computed using unit densities 
\begin{align}\label{eqn:unit_density}
e_j(x_k) = \left\{
\begin{array}{cc}
1, & k=j\\
0, & k\neq j
\end{array}
\right.
\end{align}
for $j,k=1,2,\dots,|\gamma|$. Then $j$-th column of $P_{\gamma}$ is computed as $Tr_\gamma\{G_{nh}X_{M^-}L_{nh}e_j\}$. This means we need to solve the auxiliary problem $|\gamma|$ times, which might pose a computational cost challenge in high dimensions.
\end{remark}

It is also proved that the boundary equations \eqref{eqn:bep} can be further restricted to $\gamma^+$ only, which gives the following theorem.
\begin{theorem}\label{thm:bep_reduced}
The grid function $u_\gamma$ is a trace to the solution of the discrete problem \eqref{eqn:discrete-eq}, if and only if the following boundary equations with projections hold:
\begin{align}\label{eqn:reduced-bep}
u_{\gamma^+} - P_{\gamma^+} u_\gamma = G_{nh}f_{\gamma^+}
\end{align}
where $P_{\gamma^+} u_\gamma:=Tr_{\gamma^+} P_{N^+\gamma}u_\gamma$.
\end{theorem}
The interested reader can refer to the monograph \cite{ryaben2001method}, or \cite[Section 5.3]{MR3954435} and \cite{Epshteyn_2012} for the proof of Theorem~\ref{thm:bep_reduced}.

\paragraph{The local basis function approach}
Thus far, what we have described is within the classical framework of difference potentials, and more details and theories of difference potentials can be found, for example, in the monograph \cite{ryaben2001method}. Next, we will discuss how to couple the boundary equations \eqref{eqn:bep} with the boundary conditions \eqref{eqn:bc}. We will first proceed with the 1D case \eqref{eqn:bc_1d} with second order accuracy.

Now we introduce the local basis functions defined only at point $x_j$ in $\gamma$, which are piecewise-smooth Lagrange polynomials. For example, basis functions with degree 1 polynomials are standard ``tent'' function:
\begin{equation}
\phi^{(1)}_{j}(\xi) = \begin{cases}
\dfrac{\xi+h}{h},& -h<\xi\leq0,\\[.7em]
\dfrac{\xi-h}{-h},& \ \ \: 0<\xi\leq h,\\[.7em]
\multicolumn{1}{@{}c@{\quad}}{0,} & \multicolumn{1}{@{}c@{\quad}}{\mbox{else.}} 
\end{cases}
\end{equation}
where $\xi=x-x_j$ for $x_j\in\gamma$, and the superscript $(1)$ denotes the degree of the polynomials. With second order accurate centered finite difference schemes, the number of points in $\gamma$ that are near each boundary point are exactly two, which is sufficient to construct degree 1 basis functions. On the grid points $x_j,j=1,\dots,N$, the basis function $\phi^{(1)}_{j}$ assumes a unit value at point $x_j$ and 0 elsewhere, which also coincides with the unit density for construction of the difference potentials operator $P_\gamma$. 
These basis functions are introduced in the Galerkin difference method \cite{banks2016galerkin}, and have been successful in many application \cite{banks2021discontinuous,banks2018galerkin,jacangelo2020galerkin,zhang2022energy}. 


With the basis functions defined, we next introduce the boundary cells $C_a$ and $C_b$ that contain the boundary points $x=a$ and $x=b$, respectively. We assume functions in the boundary cells are linear combinations of the basis functions, i.e.,
\begin{align}
    u(x) = \sum_{x_j\in\gamma} c_j \phi^{(1)}_j(x), \quad x\in C_a\cup C_b,
\end{align}
where $c_j$ is the unknown coefficient to be determined. Due to the unit value and interpolating property of the basis functions at point $x_j$, the coefficients $c_j$ are exactly the unknown function values $u_j$ at the discrete grid boundary $\gamma$, which appears in the boundary equation \eqref{eqn:bep}., i.e.,
\begin{align}
    c_j = u_{j}
\end{align}
for $x_j\in \gamma$. Clearly, this representation of functions in the boundary cells gives flexibility in both handling different types of boundary conditions and allowing arbitrary locations of the boundary points. For example, in the case of Dirichlet boundary conditions, we can simply enforce the boundary conditions in the collocation fashion,
\begin{align}\label{eqn:1d-O2-bc}
\sum_{x_j\in \gamma} u_j \phi^{(1)}_j(a) = u_a, \quad \sum_{x_j\in \gamma} u_j \phi^{(1)}_j(b) = u_b,
\end{align}
regardless of the position of $a$ or $b$.
Other types of boundary conditions can be treated similarly. In this work, we will use the above strong form of boundary conditions for further development. Weak formulations will be investigated in future work for more general boundary conditions.

Coupling the reduced boundary equations \eqref{eqn:reduced-bep} with the boundary condition \eqref{eqn:1d-O2-bc} leads to the following square system of equations:
\begin{subequations}\label{eqn:boundary-eqs}
\begin{align}
u_{\gamma^+}-P_{\gamma^+}u_{\gamma} &= G_hf_{\gamma^+},\\
\sum_{x_j\in \gamma} u_j \phi^{(1)}_j(a) &= u_a,\\
\sum_{x_j\in \gamma} u_j \phi^{(1)}_j(b) &= u_b.
\end{align}
\end{subequations}
from which we solve for the unknown density $u_\gamma$. After obtaining the density $u_\gamma$, the approximated solution $u_h$ is obtained using the discrete generalized Green's formula
\begin{align}\label{eqn:green}
u_h = P_{N^+\gamma} u_\gamma + G_{2h}f.
\end{align}

\begin{figure}[H]
\centering
\begin{tikzpicture}[scale=1.5]
\draw[thick,|-|] (-1.5,0) -- (0.25,0);
\draw[thick] (0.25,0) -- (3.75,0);
\draw[thick,|-|] (3.75,0) -- (5.5,0);
\filldraw (0.5,0) circle (1pt);
\filldraw (1.0,0) circle (1pt);
\filldraw (1.5,0) circle (1pt);
\filldraw (2.0,0) circle (1pt);
\filldraw (2.5,0) circle (1pt);
\filldraw (3.0,0) circle (1pt);
\filldraw (3.5,0) circle (1pt);

\draw (-1.0,0) circle (1pt);
\draw (-0.5,0) circle (1pt);
\draw (0,0) circle (1pt);

\draw (-0.5,0) +(-1.5pt,-1.5pt) rectangle +(1.5pt,1.5pt);
\draw (0,0) +(-1.5pt,-1.5pt) rectangle +(1.5pt,1.5pt);
\draw (0.5,0) +(-1.5pt,-1.5pt) rectangle +(1.5pt,1.5pt);
\draw (1,0) +(-1.5pt,-1.5pt) rectangle +(1.5pt,1.5pt);

\draw (4,0) circle (1pt);
\draw (4.5,0) circle (1pt);
\draw (5.0,0) circle (1pt);

\draw (4.5,0) +(-1.5pt,-1.5pt) rectangle +(1.5pt,1.5pt);
\draw (4,0) +(-1.5pt,-1.5pt) rectangle +(1.5pt,1.5pt);
\draw (3.5,0) +(-1.5pt,-1.5pt) rectangle +(1.5pt,1.5pt);
\draw (3,0) +(-1.5pt,-1.5pt) rectangle +(1.5pt,1.5pt);

\draw (-0.25,0.1) node[above] {$a$};
\draw (4.25,0.1) node[above] {$b$};
\draw (-1.5,0.1) node[above] {$a-\ell$};
\draw (5.5,0.1) node[above] {$b+\ell$};

\filldraw (1.0,-0.5) circle (1pt);
\draw (1.3,-0.5) node {$M^+$};
\draw (2,-0.5) circle (1pt);
\draw (2.3,-0.5) node {$M^-$};
\draw (3,-0.5) +(-1.5pt,-1.5pt) rectangle +(1.5pt,1.5pt);
\draw (3.3,-0.5) node {$\gamma$};
\end{tikzpicture}
\caption{An example of point sets with the second order accurate discretization}
\label{fig:1d-example-O4}
\end{figure}
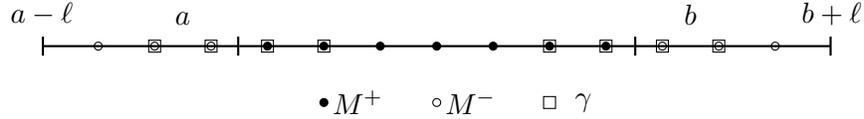

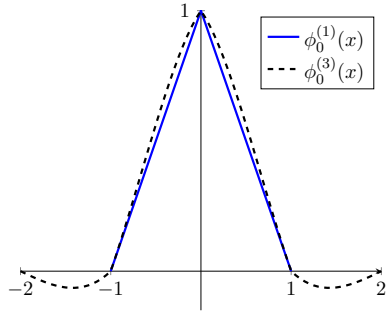
\begin{figure}[H]
    \centering
    \begin{tikzpicture}[scale=0.7,
        declare function={
        func2(\x)= (\x < -1) * (0)   +
                    and(\x >= -1, \x < 1) * (1-abs(\x))     +
                    (\x >= 1) * (0);
        func4(\x)= (\x < -1) * (0)   +
                    and(\x >= -2, \x < -1) * (\x+3)*(\x+2)*(\x+1)/6 + 
                    and(\x >= -1, \x < -0) * (\x+2)*(\x+1)*(\x-1)/(-2) + 
                    and(\x >= 0, \x < 1) * (\x+1)*(\x-1)*(\x-2)/2 + 
                    and(\x >= 1, \x < 2) * (\x-1)*(\x-2)*(\x-3)/(-6) + 
                    (\x > 2) * (0);
        }
    ]
    \begin{axis}[
        axis x line=middle, axis y line=middle,
        ymin=-0.15, ymax=1, ytick={1},
        xmin=-2, xmax=2, xtick={-2,...,2}, 
    ]
    
    \addplot[blue,very thick,samples at={-1,-0.95,...,1}, unbounded coords=jump] {func2(x)};
    \addplot[dashed,very thick,samples at={-2,-1.95,...,2}, unbounded coords=jump] {func4(x)};
    \legend{$\phi^{(1)}_0(x)$,$\phi^{(3)}_0(x)$}
    \end{axis}
    \end{tikzpicture}
    \caption{Normalized basis functions in 1D}    
    \label{fig:basis_functions}
\end{figure}

\paragraph{High order accuracy}
We will use the fourth order accurate scheme as an example to illustrate the constructions of higher order schemes. The point sets $N^0,N^\pm,M^0,M^\pm,\gamma,\gamma^\pm$ are constructed similarly as in the second order case. An example of the point sets can be found in Figure~\ref{fig:1d-example-O4}.
In the case of fourth order accuracy, we observe that there are four points (two inside and two outside) in the $\gamma$ set for each boundary point, resulting from the wider centered finite difference stencils. Thus near each boundary point, we readily obtained enough grid points to define basis functions of degree three, which are defined by
\begin{equation}
\phi^{(3)}_{j}(\xi) =
\begin{cases}
\dfrac{(\xi+3h)(\xi+2h)(\xi+h)}{6h^3},& -2h<\xi\leq-h,\\[.7em]
\dfrac{(\xi+2h)(\xi+h)(\xi-h)}{-2h^3},& \ \:-h<\xi\leq0,\\[.7em]
\dfrac{(\xi+h)(\xi-h)(\xi-2h)}{2h^3},& \quad \ \,0<\xi\leq h,\\[.7em]
\dfrac{(\xi-h)(\xi-2h)(\xi-3h)}{-6h^3},& \quad \ \,h<\xi\leq 2h,\\[.7em]
\multicolumn{1}{@{}c@{\quad}}{0,} & \multicolumn{1}{@{}c@{\quad}}{else.}
\end{cases}
\end{equation}
with $\xi=x-x_j$. A sketch of the basis functions $\phi^{(1)}_{j}$ and $\phi^{(3)}_{j}$ can be found in Figure~\ref{fig:basis_functions}.
Standard Lagrange polynomials defined at the four grid points in $\gamma$ set at each boundary point works as well in 1D, but is difficult to extend to higher dimensions. The local piecewise basis functions $\phi^{(3)}_{j}$ defined above, however, assume the same shape at each grid point in $\gamma$, which can easily to extended to higher dimensions or higher accuracy.
Another advantage of the basis functions $\phi^{(3)}_j(\xi)$ is that it is centered at $\xi=0$ and  its support is compact, i.e., the basis functions assume zero values outside of $-2h\leq\xi\leq2h$. 

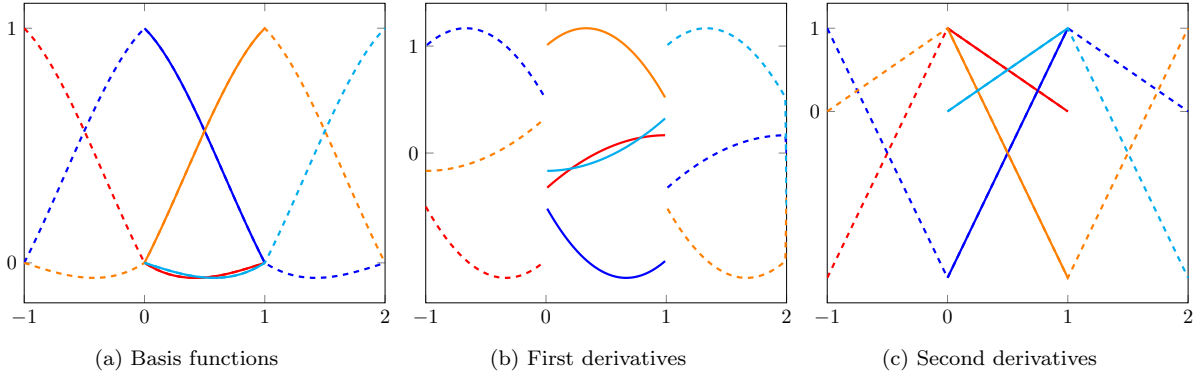
\begin{figure}
    \centering
    \begin{subfigure}{0.3\textwidth}
    \centering
    \begin{tikzpicture}[scale=0.7,
        declare function={
        func2(\x)= (\x < -1) * (0)   +
                    and(\x >= -1, \x < 1) * (1-abs(\x))     +
                    (\x >= 1) * (0);
        func4(\x,\y)= ((\x-\y) < -1) * (0)   +
                    and((\x-\y) >= -2, (\x-\y) < -1) * ((\x-\y)+3)*((\x-\y)+2)*((\x-\y)+1)/6 + 
                    and((\x-\y) >= -1, (\x-\y) < -0) * ((\x-\y)+2)*((\x-\y)+1)*((\x-\y)-1)/(-2) + 
                    and((\x-\y) >= 0, (\x-\y) < 1) * ((\x-\y)+1)*((\x-\y)-1)*((\x-\y)-2)/2 + 
                    and((\x-\y) >= 1, (\x-\y) < 2) * ((\x-\y)-1)*((\x-\y)-2)*((\x-\y)-3)/(-6) + 
                    ((\x-\y) > 2) * (0);
        }
    ]
    \begin{axis}[
        ytick={0,1},
        xmin=-1, xmax=2, xtick={-1,0,1,2}
    ]
    
    \addplot[dashed,red,very thick,samples at={-1,-0.95,...,1}, unbounded coords=jump] {func4(x,-1)};
    \addplot[red,very thick,samples at={0,0.05,...,1}, unbounded coords=jump] {func4(x,-1)};
    \addplot[dashed,blue,very thick,samples at={-1,-0.95,...,2}, unbounded coords=jump] {func4(x,0)};
    \addplot[blue,very thick,samples at={0,0.05,...,1}, unbounded coords=jump] {func4(x,0)};
    \addplot[dashed,orange,very thick,samples at={-1,-0.95,...,2}, unbounded coords=jump] {func4(x,1)};
    \addplot[orange,very thick,samples at={0,0.05,...,1}, unbounded coords=jump] {func4(x,1)};
    \addplot[dashed,cyan,very thick,samples at={0,0.05,...,2}, unbounded coords=jump] {func4(x,2)};
    \addplot[cyan,very thick,samples at={0,0.05,...,1}, unbounded coords=jump] {func4(x,2)};
    \end{axis}
    \end{tikzpicture}
    \caption{Basis functions}\label{fig:boundary_cell}
    \end{subfigure}
    ~
    \begin{subfigure}{0.3\textwidth}
    \centering
    \begin{tikzpicture}[scale=0.7,
        declare function={
        func2(\x)= (\x < -1) * (0)   +
                   (\x < 0 ) * (1)
                    and(\x < 1) * (-1)     +
                    (\x >= 1) * (0);
        func4(\x,\y)= ((\x-\y) < -1) * (0)   +
                    and((\x-\y) >= -2, (\x-\y) < -1) * ( ((\x-\y)+2)*((\x-\y)+1)/6 + ((\x-\y)+3)*((\x-\y)+1)/6 + ((\x-\y)+3)*((\x-\y)+2)/6) + 
                    and((\x-\y) >= -1, (\x-\y) < -0) * ( ((\x-\y)+1)*((\x-\y)-1)/(-2) + ((\x-\y)+2)*((\x-\y)-1)/(-2) + ((\x-\y)+2)*((\x-\y)+1)/(-2) ) + 
                    and((\x-\y) >= 0, (\x-\y) < 1) * ( ((\x-\y)-1)*((\x-\y)-2)/2 + ((\x-\y)+1)*((\x-\y)-2)/2 + ((\x-\y)+1)*((\x-\y)-1)/2 ) + 
                    and((\x-\y) >= 1, (\x-\y) < 2) * ( ((\x-\y)-2)*((\x-\y)-3)/(-6) + ((\x-\y)-1)*((\x-\y)-3)/(-6)+ ((\x-\y)-1)*((\x-\y)-2)/(-6) ) + 
                    ((\x-\y) > 2) * (0);
        }
    ]
    \begin{axis}[
        ytick={0,1},
        xmin=-1, xmax=2, xtick={-1,0,1,2}
    ]
    
    \addplot[dashed,red,very thick,samples at={-1,-0.99,...,-0.01}, unbounded coords=jump] {func4(x,-1)};
    \addplot[red,very thick,samples at={0.01,0.02,...,0.99}, unbounded coords=jump] {func4(x,-1)};
    \addplot[dashed,blue,very thick,samples at={-0.99,-0.98,...,-0.01}, unbounded coords=jump] {func4(x,0)};
    \addplot[blue,very thick,samples at={0.01,0.02,...,0.99}, unbounded coords=jump] {func4(x,0)};
    \addplot[dashed,blue,very thick,samples at={1.01,1.02,...,2}, unbounded coords=jump] {func4(x,0)};
    \addplot[dashed,orange,very thick,samples at={-1,-0.99,...,-0.01}, unbounded coords=jump] {func4(x,1)};
    \addplot[orange,very thick,samples at={0.01,0.02,...,0.99}, unbounded coords=jump] {func4(x,1)};
    \addplot[dashed,orange,very thick,samples at={1.01,1.02,...,2}, unbounded coords=jump] {func4(x,1)};
    \addplot[dashed,cyan,very thick,samples at={1.01,1.02,...,2}, unbounded coords=jump] {func4(x,2)};
    \addplot[cyan,very thick,samples at={0.01,0.02,...,0.99}, unbounded coords=jump] {func4(x,2)};
    \end{axis}
    \end{tikzpicture}
    \caption{First derivatives}\label{fig:1st_deriv}
    \end{subfigure}
    ~
    \begin{subfigure}{0.3\textwidth}
    \centering
    \begin{tikzpicture}[scale=0.7,
        declare function={
        func2(\x)= (\x < -1) * (0)   +
                    and(\x >= -1, \x < 1) * (0)     +
                    (\x >= 1) * (0);
        func4(\x,\y)= ((\x-\y) < -1) * (0)   +
                    and((\x-\y) >= -2, (\x-\y) < -1) * ( 2*((\x-\y)+3)/6 + 2*((\x-\y)+2)/6 + 2*((\x-\y)+1)/6 ) + 
                    and((\x-\y) >= -1, (\x-\y) < -0) * ( 2*((\x-\y)+2)/(-2) + 2*((\x-\y)+1)/(-2) + 2*((\x-\y)-1)/(-2) ) + 
                    and((\x-\y) >= 0, (\x-\y) < 1) * ( ((\x-\y)+1) + ((\x-\y)-1) + ((\x-\y)-2) ) + 
                    and((\x-\y) >= 1, (\x-\y) < 2) * ( 2*((\x-\y)-1)/(-6) + 2*((\x-\y)-2)/(-6) + 2*((\x-\y)-3)/(-6) )+ 
                    ((\x-\y) > 2) * (0);
        }
    ]
    \begin{axis}[
        ytick={0,1},
        xmin=-1, xmax=2, xtick={-1,0,1,2}
    ]
    
    \addplot[dashed,red,very thick,samples at={-1,-0.95,...,1}, unbounded coords=jump] {func4(x,-1)};
    \addplot[red,very thick,samples at={0,0.05,...,1}, unbounded coords=jump] {func4(x,-1)};
    \addplot[dashed,blue,very thick,samples at={-1,-0.95,...,2}, unbounded coords=jump] {func4(x,0)};
    \addplot[blue,very thick,samples at={0,0.05,...,1}, unbounded coords=jump] {func4(x,0)};
    \addplot[dashed,orange,very thick,samples at={-1,-0.95,...,2}, unbounded coords=jump] {func4(x,1)};
    \addplot[orange,very thick,samples at={0,0.05,...,1}, unbounded coords=jump] {func4(x,1)};
    \addplot[dashed,cyan,very thick,samples at={0,0.05,...,2}, unbounded coords=jump] {func4(x,2)};
    \addplot[cyan,very thick,samples at={0,0.05,...,1}, unbounded coords=jump] {func4(x,2)};
    \end{axis}
    \end{tikzpicture}
    \caption{Second derivatives}\label{fig:2nd_deriv}
    \end{subfigure}
    \caption{Basis functions and their derivatives}\label{fig:basis_and_deriv}
\end{figure}

In Figure~\ref{fig:basis_and_deriv}, sketches of all the basis functions defined at these four points and their first and second derivatives are given. We assume points $-1,0,1,2$ are the four grid points in the $\gamma$ set and the boundary point $a$ or $b$ falls into the middle cell $[0,1]$. The boundary cell $[0,1]$ is exactly spanned by Lagrange basis functions of degree three, which is what we need, though the other two cells are are not. Boundary closures using these basis functions defined at $\gamma$ points fit exactly into the difference potentials framework. 

Similarly as in the second order case, we obtain the following boundary equations:
\begin{subequations}\label{eqn:boundary-eqs-O4}
\begin{align}
u_{\gamma^+}-P_{\gamma^+}u_{\gamma} = G_{4h}f_{\gamma^+},\\
\sum_{x_j\in\gamma} u_j \phi^{(3)}_j(a) = u_a,\\
\sum_{x_j\in\gamma} u_j \phi^{(3)}_j(b) = u_b.
\end{align}
\end{subequations}
The operators $P_{\gamma^+}$ and $G_{4h}$ now are formulated using fourth order accurate finite difference schemes. In detail, the $j$-th column of operator $P_{\gamma^+}$ is computed as $Tr_{\gamma}G_{4h}X_{M^-}L_{4h}e_j$ where $L_{4h}$ denotes also the fourth order accurate discrete version.
However, the equations in \eqref{eqn:boundary-eqs-O4} are under-determined, and we need to supply additional two more boundary equations in order to balance the number of unknowns and the number of equations. Compatibility condition is one solution, i.e., the following two more equations are added to \eqref{eqn:boundary-eqs-O4} so as to obtain a square system of equations for Poisson's equation:
\begin{subequations}\label{eqn:boundary-eqs-compat}
\begin{align}
\sum_{x_j\in \gamma} u_j \frac{d^2\phi^{(3)}_j}{dx^2}(a) = f(a),\\
\sum_{x_j\in \gamma} u_j \frac{d^2\phi^{(3)}_j}{dx^2}(b) = f(b).
\end{align}
\end{subequations}
Then we solve \eqref{eqn:boundary-eqs-O4} and \eqref{eqn:boundary-eqs-compat} for the density $u_\gamma$, and the approximated solution is obtained using the discrete generalized Green's formula \eqref{eqn:green} with fourth order accurate operators.

Higher order accurate versions will follow the fourth order case. The major difference will be the width of grid points in $\gamma$ set and the degree of piecewise Lagrange polynomials, plus a high order compatibility condition.

\subsection{Two dimensions}

In this part, we consider a 2D domain $\Omega$ of arbitrary shape, and Dirichlet boundary condition $u(x)=g(x)$ is assumed on the boundary $\partial\Omega$.

\paragraph{Second order accuracy}
Again we start with the introduction of an auxiliary domain and discretization of the auxiliary domain with uniform Cartesian meshes. As in the 1D case, we will introduce the point sets $M^0,M^\pm,N^0,N^\pm$ and the discrete grid boundary $\gamma$ (and $\gamma^\pm$) based on the 2D centered finite difference stencils
\begin{align}\label{eqn:2d_stencil}
\mathcal{N}_{j,k}^{n} = \bigcup_{l_1=0}^n\bigcup_{l_2=0}^n x_{j\pm l_1,k\pm l_2}
\end{align} 
with $n=1$ for second order accuracy and $n=2$ for fourth order accuracy.

With the uniform grids defined, we discretize the PDE in $M^+$ and reduce the discretization to boundary equations with projections in $\gamma^+$ using Theorem~\ref{thm:bep_reduced}. Then we obtain the following boundary equations for the second order discretization as in 1D,
\begin{subequations}\label{eqn:boundary-eqs-2d-2nd}
\begin{align}
u_{\gamma^+}-P_{\gamma^+}u_{\gamma} = G_{2h}f_{\gamma^+}\label{eqn:boundary-eqs-2d-2nd-1}\\
\sum_{x_{j,k}\in\gamma} u_{j,k}\phi^{(1,1)}_{j,k}(x_b) = g(x_b)\label{eqn:boundary-eqs-2d-2nd-2}
\end{align}
\end{subequations}
The difference potentials operator $P_{\gamma^+}$ is constructed similarly as in 1D, namely the $j$-th column of $P_{\gamma^+}$ is computed as $Tr_{\gamma^+}G_{2h}X_{M^-}L_{2h}e_j$ where $e_j$ is the unit density \eqref{eqn:unit_density}, except that the exact definitions of $e_j$, $G_{2h}$ and $L_{2h}$ have changed due to the dimension change. The local basis function $\phi^{(1,1)}_{j,k}$ are defined at each point $x_{j,k}\in\gamma$, and is obtained as the tensor products of 1D local basis functions,
\begin{align}\label{eqn:basis_function_2d}
\phi^{(1,1)}_{j,k}=\phi^{(1)}_{j}\phi^{(1)}_k,
\end{align}
where $\phi^{(1)}_j$ and $\phi^{(1)}_{k}$ are the 1D local linear basis function in the $x$ and $y$ directions respectively.

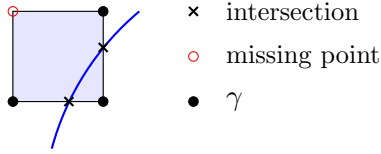
\begin{figure}[H]
    \centering
    \begin{tikzpicture}[scale=1.2]
    \draw[step=1cm,gray,very thin] (0,0) grid (1,1);
    \filldraw[fill=blue!10] (0,0) rectangle (1,1);
    \draw[blue,thick] (1.4,1) arc (130:165:3);
    \filldraw (0,0) circle (1.5pt);
    \filldraw (1,0) circle (1.5pt);
    \draw[red] (0,1) circle (1.5pt);
    \filldraw (1,1) circle (1.5pt);
    \draw[mark size=+2pt,thick] plot[mark=x] coordinates {(0.62,0)};
    \draw[mark size=+2pt,thick] plot[mark=x] coordinates {(1,0.6)};
    \draw[mark size=+2pt,thick] plot[mark=x] coordinates {(2,1)};
    \draw (2.25,1) node[right] {\small intersection};
    \filldraw (2,0) circle (1.5pt);
    \draw (2.25,0) node[right] {$\gamma$};
    \draw[red] (2,0.5) circle (1.5pt);
    \draw (2.25,0.5) node[right] {\small missing point};
    \end{tikzpicture}
    \caption{Missing point and intersection points}
    \label{fig:missing_point}
    \end{figure}

The point $x_b$ in \eqref{eqn:boundary-eqs-2d-2nd-2} denotes some boundary point on $\partial\Omega$. In 1D, it was straightforward to impose the boundary condition at the two boundary points. Here in 2D, we are faced with multiple choices of the location of the boundary points. A key observation in the 1D case was that the function space in the boundary cells can be spanned by a complete set of the Lagrange basis functions, while other cells do not admit this property. Now in 2D, as can be seen in Figure~\ref{fig:missing_point}, the shaded boundary cell is missing the red point. In the second order case, one natural choice to mitigate the effect of missing points is to select $x_b$ as the intersection point between the grid lines and the boundary curve, since the contribution from the basis function at the missing point vanishes at the two further grid lines that include the intersection points, i.e., $\phi_j(x_b)=0$. 

\begin{figure}[H]
\centering
\begin{tikzpicture}[scale=1.2]
\draw[step=1cm,gray,very thin] (0,0) grid (1,1);
\draw[blue,thick] (0.6,1.2) arc (140:165:4);
\filldraw (0,1) circle (1.5pt);
\draw (0.25,-1.5) node[right] {intersection};
\draw[mark size=+2pt,thick] plot[mark=x] coordinates {(0,-1.5)};
\draw (0.5,-0.5) node[below] {1-to-2};
\draw[mark size=+2pt,thick] plot[mark=x] coordinates {(0,0.25)};
\draw[mark size=+2pt,thick] plot[mark=x] coordinates {(0.45,1)};
\draw[step=1cm,gray,very thin] (3,0) grid (4,1);
\draw[blue,thick] (0.9+3,1.2) arc (140:165:4);
\filldraw (3,1) circle (1.5pt);
\draw[mark size=+2pt,thick] plot[mark=x] coordinates {(3.75,1)};
\filldraw (3,-1.5) circle (1.5pt);
\draw (3.25,-1.5) node[right] {$\gamma^-$};
\draw (0.5+3,-0.5) node[below] {1-to-1};
\end{tikzpicture}
\caption{$\gamma^-$ vs intersection}
\label{fig:intersection_O2}
\end{figure}
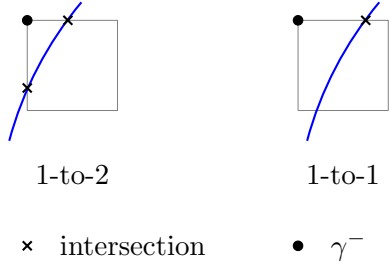

There are several advantages of using the intersection points. First, explicit or implicit geometry can be treated in the same fashion, as long as the intersection points can be found. Other choices, such as equi-distant points, might require explicit representation of the geometry, or difficult to extend to more complicated geometry, especially in 3D. Second, to find the intersection points, one only need to solve multiple 1D root finding problems, regardless of the dimensions, thus can be easily extended to 3D. If every such intersection points are included for \eqref{eqn:boundary-eqs-2d-2nd-2}, the boundary equations \eqref{eqn:boundary-eqs-2d-2nd} would be a least square system, which might be more expensive to solve. On the other hand, square systems might be cheaper and easier to solve. To balance the number of unknowns and the number of equations, we first note that the total number of unknowns in \eqref{eqn:boundary-eqs-2d-2nd} again is $|\gamma|$. The first equation \eqref{eqn:boundary-eqs-2d-2nd-1} contributes $|\gamma^+|$ number of equations, so we would expect $|\gamma^-|=|\gamma|-|\gamma^+|$ number of equations in the second equation \eqref{eqn:boundary-eqs-2d-2nd-2}, which means we need $|\gamma^-|$ boundary points to impose the boundary conditions. By noting the 1-to-2 or 1-to-1 relation (Figure \ref{fig:intersection_O2}) of points in $\gamma^-$ with all the boundary points, we select an intersecting point for each point in $\gamma^-$ to ensure the 1-to-1 correspondence. For the 1-to-2 case, either intersection point works and we will use the nearest one. Basically, we can regard \eqref{eqn:boundary-eqs-2d-2nd-1} as accounting for the point set $\gamma^+$, and \eqref{eqn:boundary-eqs-2d-2nd-2} for the point set $\gamma^-$.

However, the approach of using only basis functions defined at $\gamma$ set and the system of boundary equations \eqref{eqn:boundary-eqs-2d-2nd} with boundary conditions imposed at the boundary intersecting points, will only work for Dirichlet boundary condition, as the derivatives of the basis functions will not vanish at the boundary, and the missing points lead to an incomplete set of basis functions to resolve the boundary conditions. In addition, the system \eqref{eqn:boundary-eqs-2d-2nd} will not be able to be extended to high order accuracy due to the same issue, i.e. more missing points lead to incomplete span of basis functions.

\begin{figure}[H]
\centering
\begin{tikzpicture}[scale=1.2]
\draw[step=1cm,gray,very thin] (0,0) grid (4,4);
\filldraw[fill=blue!10] (1,2) rectangle (2,3);
\draw (1.5,2.6) node {$C_i$};
\draw[blue,thick] (0,0.2) .. controls (1,3) and (2,2) .. (4,3.8);
\filldraw (0,1) circle (1.5pt);
\filldraw (1,1) circle (1.5pt);
\filldraw (1,2) circle (1.5pt);
\filldraw (2,2) circle (1.5pt);
\filldraw (2,3) circle (1.5pt);
\draw[mark size=+2pt,thick] plot[mark=x] coordinates {(3,2)};
\draw[mark size=+2pt,thick] plot[mark=x] coordinates {(2,1)};
\draw[mark size=+2pt,thick] plot[mark=x] coordinates {(1,0)};
\filldraw (3,3) circle (1.5pt);
\filldraw (3,4) circle (1.5pt);
\filldraw (0,0) circle (1.5pt);
\filldraw (4,3) circle (1.5pt);
\filldraw (4,4) circle (1.5pt);
\draw[red] (0,2) circle (1.5pt);
\draw[red] (1,3) circle (1.5pt);
\draw[red] (2,4) circle (1.5pt);
\draw[red] (5,3) circle (1.5pt);
\draw[red] (5.25,3) node[right] {$\eta$};
\filldraw (5,2) circle (1.5pt);
\filldraw (5.25,2) node[right] {$\gamma$};
\draw[mark size=+2pt,thick] plot[mark=x] coordinates {(5,1)};
\draw (5.25,1) node[right] {$\omega$};
\end{tikzpicture}
\caption{Example of $\eta,\gamma,\omega$ and a cut cell $C_i$ }
\label{fig:5-point-example}
\end{figure}
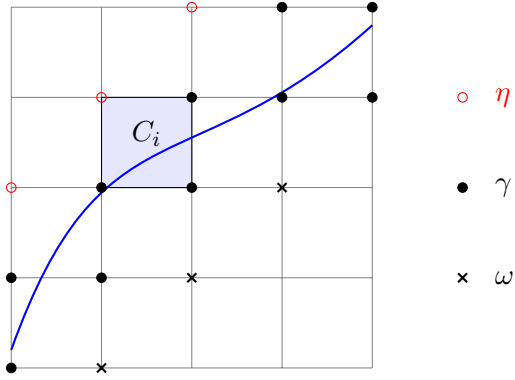

To this end, we introduce two more point sets that complement the discrete grid boundary $\gamma$.
\begin{definition} \label{def:eta_omega_O2}
The minimal grid point sets together with $\gamma$ that span the  $\cup_iC_i$ in second order accuracy with bilinear basis functions $\phi^{(1,1)}_{j,k}$ are:
\begin{subequations}
\begin{align}
    \eta:=M^-\cap\left(\bigcup_{i}V(C_i) \right)\backslash \:\gamma^-,\\
    \omega:=M^+\cap\left(\bigcup_{i}V(C_i) \right)\backslash \: \gamma^+,
\end{align}
\end{subequations}
where $C_i$ denotes a boundary cell that intersects with the boundary.
\end{definition}

\begin{remark}
    In practice, we do not need to include every boundary cell that intersects with the boundary $\partial\Omega$. Instead, we will only need those cells that are support of the intersection points that will go into the boundary system \eqref{eqn:boundary-eqs-2d-2nd}. To do this, we note that the piecewise smooth basis functions include the right end points but not the left ones. Thus, for an intersection point $x_b$, we will use the boundary cell to the left or to the below, which is indicated in the figure below.
    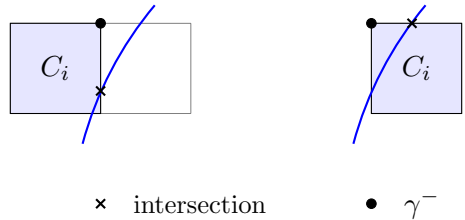
\begin{figure}[H]
    \centering
    \begin{tikzpicture}[scale=1.2]
    \draw[step=1cm,gray,very thin] (-1,0) grid (1,1);
    \filldraw[fill=blue!10] (-1,0) rectangle (0,1);
    \draw[blue,thick] (0.6,1.2) arc (140:165:4);
    \filldraw (0,1) circle (1.5pt);
    \draw (-0.5,0.5) node {$C_i$};
    \draw[mark size=+2pt,thick] plot[mark=x] coordinates {(0,0.25)};
    \draw[step=1cm,gray,very thin] (3,0) grid (4,1);
    \filldraw[fill=blue!10] (3,0) rectangle (4,1);
    \draw[blue,thick] (0.6+3,1.2) arc (140:165:4);
    \filldraw (3,1) circle (1.5pt);
    \draw (3.5,0.5) node {$C_i$};
    \draw[mark size=+2pt,thick] plot[mark=x] coordinates {(3.45,1)};
    \draw[mark size=+2pt,thick] plot[mark=x] coordinates {(0,-1)};
    \draw (0.25,-1) node[right] {\small intersection};
    \filldraw (3,-1) circle (1.5pt);
    \draw (3.25,-1) node[right] {$\gamma^-$};
    \end{tikzpicture}
    \caption{Intersection points and their support cells}
    \label{fig:support_cell}
    \end{figure}
    Figure~\ref{fig:support_cell} shows that the choices of intersection points will determine the supporting cell, thus the point set $\eta$ and $\omega$ from the vertices of the supporting cell.
\end{remark}

With the discrete grid boundary augmented by $\eta$ and $\omega$, we would need to reformulate the boundary equations \eqref{eqn:boundary-eqs-2d-2nd} accordingly. To start with, we need to see the effects of adding more points to the discrete grid boundary.
\begin{proposition}\label{prop:new_bep}
BEP~\eqref{eqn:reduced-bep} is equivalent to  in the augmented discrete grid boundary $\zeta:=\gamma\cup\eta\cup\omega$
\begin{align}\label{eqn:new_bep}
    u_{\zeta^+}-P_{\zeta^+}u_{\zeta} = G_{2h}f_{\zeta^+}.
\end{align}
\end{proposition}
\begin{proof}
By the definition, for any density $u_\omega$ defined on $\omega$, its difference potentials is given by
\begin{align}
P_{N^+\omega}u_{\omega} = G_{2h}X_{M^-}L_{2h}u_\omega = 0,
\end{align}
since the stencil defined in $\omega$ does not reach out of $M^+$, which implies $P_{\zeta^+}u_\omega=0$. Then for any density $u_\eta$ defined on set $\eta$, we have
\begin{align}
P_{N^+\eta}u_{\eta} = G_{2h}X_{M^-}L_{2h}u_\eta = G_{2h}L_{2h}u_\eta = u_\eta
\end{align}
since the stencils in $\eta$ do not enter $M^+$, which implies $P_{\zeta^+}u_\eta=0$ on $\zeta^+$ since $u_\eta$ can only be nonzero on $\eta$ set. 

To conclude, we are only adding trivial information $u_{\omega}=G_{2h}f_\omega$ to the BEP~\eqref{eqn:reduced-bep}, without changing its structure. 
\end{proof}

Proposition~\ref{prop:new_bep} states that the addition of point set $\omega$ does not affect the structure of boundary equations, and is more of a notation convenience to include the extra points in $\eta$ and $\omega$. On the other hand, the ghost values at $\eta$ can be extrapolated from the density $u_\gamma$ in $\gamma$, plus $u_\omega$ if higher order extrapolation is needed. We will denote the extrapolation as:
\begin{align}
    Ru_{\eta\cup\gamma\cup\omega} = 0
\end{align}
where the extrapolation operator $R$ is be constructed, for example, from the second order extrapolation $u_{-1}=2u_0-u_1$ in $x$ or $y$ direction. A more robust choice is to employ averaged extrapolations from both $x$ and $y$ directions.
Now the boundary conditions need to be reformulated with all basis functions defined on the augmented discrete grid boundary $\zeta$:
\begin{align}
    \sum_{x_{j,k}\in \zeta} u_{j,k}\phi^{(1,1)}_{j,k}(x_b) = g(x_b)
\end{align}
where $x_b$ denotes some intersecting point.

Now we are ready to assemble the effective boundary equations:
\begin{align}\label{eqn:boundary_eq_O2}
    \left(
    \begin{array}{c c}
        0 & I-P_{\zeta^{+}} \\
        \multicolumn{2}{c}{\Phi^{(1,1)}(x_b)} \\
        \multicolumn{2}{c}{R}\\
    \end{array}
    \right)
    \begin{pmatrix}
        u_\eta\\
        u_\gamma\\
        u_\omega
    \end{pmatrix}
=
    \begin{pmatrix}
        G_hf_{\zeta^{+}}\\
        g(x_b)\\
        0
    \end{pmatrix}
\end{align}
where the first row block corresponds to \eqref{eqn:new_bep}, the second row for the boundary condition~\eqref{eqn:boundary-eqs-2d-2nd-2} with $\Phi^{(1,1)}$ assembled from the basis functions and the last row from the extrapolations. The operator $I$ is not a square matrix, but rather a part of the identity matrix with the same size of $P_{\zeta^+}$.

The total number of unknowns in the boundary equations \eqref{eqn:boundary_eq_O2} are $|\zeta|=|\eta|+|\gamma|+|\omega|$. The first row contributes rank $|\gamma^+|+|\omega|$ and the third row $|\eta|$, which leaves $|\gamma^-|$ number of points on the boundary to impose the boundary condition for the second row. As discussed above, for each point $x_{jk}$ in $\gamma^-$, we select an intersecting point between the boundary and the grid lines that $x_j$ lies on to impose the boundary conditions, which will provide exactly $|\gamma^-|$ number of equations for the second row.

\begin{remark}
For other types of boundary conditions, the second row block will be modified accordingly
\begin{align}
\sum_{x_{j,k}\in \zeta} u_{jk}\ell (\phi^{(1,1)}_{j,k})(x_b) = g(x_b)
\end{align}
where $\ell$ denotes a boundary operator. A caveat here is that for a general boundary condition, collocation methods will result in a loss of order in accuracy, so this is not the optimal approach. The optimal treatment of the general boundary conditions is yet to be investigated. Besides, as can be observed in Figure~\ref{fig:1st_deriv}, the first derivatives of the local basis functions are not continuous at the grid lines, which means the point $x_b$ where we impose the boundary condition needs to nudge into the the boundary cell along the boundary curve for the basis functions to resolve the boundary conditions.
\end{remark}

Lastly we solve the boundary equations \eqref{eqn:boundary_eq_O2} for the density $u_\gamma$, and the approximated solution is obtained using the discrete generalized Green's formula \eqref{eqn:green}, with modifications of operators of two dimensional second order accuracy.

\paragraph{Fourth order accuracy}

For the fourth order accurate schemes, we will only comment on the difference from the second order accurate ones. The structure of point sets $M^\pm,N^\pm,\gamma$ or $\gamma^\pm$ will be similar, except a wider 9-point stencil $\mathcal{N}_{jk}^2$ is employed. Basis functions of degree three will also be constructed as tensor products:
\begin{align}
\phi^{(3,3)}_{j,k}=\phi^{(3)}_{j}\phi^{(3)}_k.
\end{align}
\begin{figure}[H]
\centering
\begin{tikzpicture}[scale=1]
\filldraw[fill=red!5] (0,1) rectangle (3,4);
\draw[step=1cm,gray,very thin] (0,0) grid (4,4);
\filldraw[fill=blue!10] (1,2) rectangle (2,3);
\draw (1.5,2.6) node {$C_i$};
\draw[blue,thick] (0,0.2) .. controls (1,3) and (2.8,1.2) .. (4,3.8);
\filldraw (0,1) circle (1.5pt);
\filldraw (1,1) circle (1.5pt);
\filldraw (2,1) circle (1.5pt);
\filldraw (0,2) circle (1.5pt);
\filldraw (1,2) circle (1.5pt);
\filldraw (2,2) circle (1.5pt);
\filldraw (3,2) circle (1.5pt);
\filldraw (1,3) circle (1.5pt);
\filldraw (2,3) circle (1.5pt);
\filldraw (3,3) circle (1.5pt);
\filldraw (4,3) circle (1.5pt);
\filldraw (3,4) circle (1.5pt);
\filldraw (1,0) circle (1.5pt);
\filldraw (2,4) circle (1.5pt);
\filldraw (0,0) circle (1.5pt);
\filldraw (4,4) circle (1.5pt);
\filldraw (3,1) circle (1.5pt);
\filldraw (4,2) circle (1.5pt);
\draw[mark size=+2pt,thick] plot[mark=x] coordinates {(3,0)};
\draw[mark size=+2pt,thick] plot[mark=x] coordinates {(2,0)};
\draw[mark size=+2pt,thick] plot[mark=x] coordinates {(4,1)};
\draw[red] (1,4) circle (1.5pt);
\draw[red] (0,3) circle (1.5pt);
\draw[red] (0,4) circle (1.5pt);
\draw[red] (5,3) circle (1.5pt);
\draw[red] (5.25,3) node[right] {$\eta$};
\filldraw (5,2) circle (1.5pt);
\filldraw (5.25,2) node[right] {$\gamma$};
\draw[mark size=+2pt,thick] plot[mark=x] coordinates {(5,1)};
\draw (5.25,1) node[right] {$\omega$};
\end{tikzpicture}
\caption{Example of $\eta,\gamma,\omega$ and a cut cell $C_i$ in fourth order accuracy}
\label{fig:9-point-example}
\end{figure}
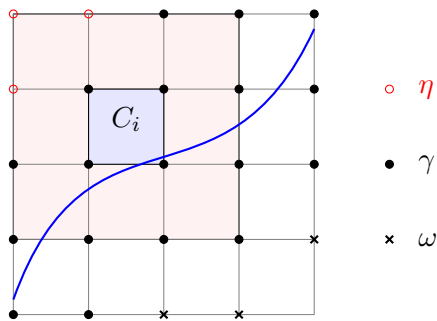

Examples of the point sets $\gamma$, $\eta$ and $\omega$ are illustrated in Figure~\ref{fig:9-point-example}. The complementary point sets $\eta$ and $\omega$ will also be defined similarly as in Definition~\ref{def:eta_omega_O2}, which are needed to establish a complete set of basis functions for the boundary cell $C_i$.
\begin{definition} \label{def:eta_omega_O4}
The minimal grid point sets together with $\gamma$ that span the cut cells $\cup_iC_i$ in fourth order accuracy using bicubic basis functions $\phi^{(3,3)}_{j,k}$ are given by
\begin{subequations}
\begin{align}
    \eta:=\left(\bigcup_{i}\bigcup_{j,k=-1,0,1}V(C_i)+(jh,kh) \right)\backslash \; \gamma \cap M^-,\\
    \omega:=\left(\bigcup_{i}\bigcup_{j,k=-1,0,1}V(C_i)+(jh,kh) \right)\backslash \; \gamma \cap M^+,
\end{align}
\end{subequations}
where $C_i$ denotes a boundary cell that intersects with the boundary.
\end{definition}
With the definition of $\eta$ and $\omega$ in Definition~\ref{def:eta_omega_O4}, we obtain a complete span for each boundary cell $C_i$. We would also like to comment on the difference of the point set $\eta$ between the second order accurate case (O2) and the fourth order accurate one (O4). In O2, function values in $\eta$ can be explicitly extrapolated using grid functions defined on $\gamma$ and $\omega$. While in O4, some of the points can only be implicitly extrapolated. Again, it can be noted that all points in $\eta$ can be extrapolated using at least 4 points in certain directions in the $\gamma$ set. We will reconcile the extrapolation between $\eta$ and $\gamma$ (and possibly $\omega$) using the extrapolation operator $R$:
\begin{align}
    Ru_{\eta\cup\gamma\cup\omega}=0,
\end{align}
The fourth order extrapolation in one direction is given by
\begin{align}
u_{-1} = -4u_0 +6 u_1-4u_2+u_3,
\end{align}
where $u_{-1}$ denotes the ghost value and the rest interior values. 
Note that the extrapolation only needs four interior points, which is exactly the width of the $\gamma$ set in the O4 case when the 9-point stencil is employed.

The unknown density now are defined on $\zeta:=\eta\cup\gamma\cup\omega$. Function values $\eta$ extrapolated from $u_{\gamma}$ (and possibly $\omega$). The set $\omega$ is incorporated into the boundary equations as stated in Proposition~\ref{prop:new_bep}. What is left is how and where to impose the boundary conditions. As in the O2 case, Dirichlet boundary condition can be imposed in the collocation fashion:
\begin{align}\label{eqn:bc1_O4}
    \sum_{x_{j,k}\in \zeta} u_{jk}\phi^{(3,3)}_{j,k}(x_b) = g(x_b),
\end{align}
where $x_b$ is some intersection point between the boundary and the grid lines, and the basis functions are again obtained using tensor products $\phi^{(3,3)}_{j,k}=\phi^{(3)}_j\phi^{(3)}_k$.

\begin{figure}[H]
\centering
\begin{tikzpicture}[scale=1]
\draw[step=1cm,gray,very thin] (0,0) grid (4,4);
\draw[blue,thick] (0,0.2) .. controls (1,3) and (2.8,1.2) .. (4,3.8);
\filldraw (0,2) circle (1.5pt);
\filldraw (1,3) circle (1.5pt);
\filldraw (2,4) circle (1.5pt);
\filldraw (3,4) circle (1.5pt);
\draw[mark size=+2pt,thick] plot[mark=asterisk] coordinates {(0,1)};
\draw[mark size=+2pt,thick] plot[mark=asterisk] coordinates {(1,2)};
\draw[mark size=+2pt,thick] plot[mark=asterisk] coordinates {(2,3)};
\draw[mark size=+2pt,thick] plot[mark=asterisk] coordinates {(3,3)};
\draw[mark size=+2pt,thick] plot[mark=asterisk] coordinates {(4,4)};
\draw[mark size=+2pt,thick] plot[mark=x] coordinates {(0.37,1)};
\draw[mark size=+2pt,thick] plot[mark=x] coordinates {(1,1.68)};
\draw[mark size=+2pt,thick] plot[mark=x] coordinates {(2,2.1)};
\draw[mark size=+2pt,thick] plot[mark=x] coordinates {(3.52,3)};
\draw[mark size=+2pt,thick] plot[mark=x] coordinates {(4,3.8)};
\draw[red,thick] (1,1.68) circle (4pt);
\draw[red,thick] (2,2.1) circle (4pt);
\draw[step=1cm,gray,very thin] (0+5,0) grid (4+5,4);
\draw[blue,thick] (0+5,0.2) .. controls (1+5,3) and (2.8+5,1.2) .. (4+5,3.8);
\filldraw (0+5,2) circle (1.5pt);
\filldraw (1+5,3) circle (1.5pt);
\filldraw (2+5,4) circle (1.5pt);
\filldraw (3+5,4) circle (1.5pt);
\draw[mark size=+2pt,thick] plot[mark=asterisk] coordinates {(0+5,1)};
\draw[mark size=+2pt,thick] plot[mark=asterisk] coordinates {(1+5,2)};
\draw[mark size=+2pt,thick] plot[mark=asterisk] coordinates {(2+5,3)};
\draw[mark size=+2pt,thick] plot[mark=asterisk] coordinates {(3+5,3)};
\draw[mark size=+2pt,thick] plot[mark=asterisk] coordinates {(4+5,4)};
\draw[mark size=+2pt,thick] plot[mark=x] coordinates {(1.7+5,2)};
\draw[mark size=+2pt,thick] plot[mark=x] coordinates {(1+5,1.68)};
\draw[mark size=+2pt,thick] plot[mark=x] coordinates {(2+5,2.1)};
\draw[mark size=+2pt,thick] plot[mark=x] coordinates {(3+5,2.525)};
\draw[red,thick] (1+5,1.68) circle (4pt);
\draw[red,thick] (2+5,2.1) circle (4pt);
\draw[mark size=+2pt,thick] plot[mark=x] coordinates {(0+1,-1)};
\draw (0.25+1,-1) node[right] {\small intersection};
\filldraw (3+1,-1) circle (1.5pt);
\draw (3.25+1,-1) node[right] {$\gamma^-_2$};
\draw[mark size=+2pt,thick] plot[mark=asterisk] coordinates {(5+1,-1)};
\draw (5.25+1,-1) node[right] {$\gamma^-_1$};
\end{tikzpicture}
\caption{Layers of $\gamma^-$ and intersection points for first layer $\gamma^-_1$ (left) and second layer $\gamma^{-}_2$ (right); shared points are highlighted in big red circle}
\label{fig:layers_gamma}
\end{figure}
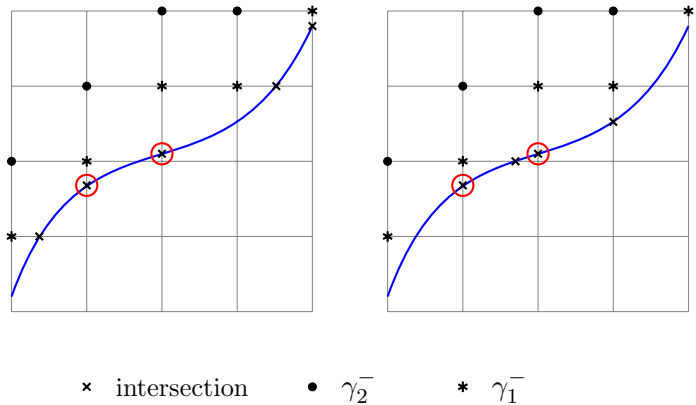

Next, we need to decide how many intersecting points to impose the BC. First, we identify $\gamma^-$ in O4 into two layers: the first layer will be identical to the discrete grid boundary in O2, which we shall denote as $\gamma^{-}_1$ and the second layer as the difference $\gamma^{-}_2:=\gamma^-\backslash\gamma^{-}_1$. For $\gamma^{-}_1$, we select an intersecting point $x_b$ as in O2 and for $\gamma^{-}_2$, we also choose an intersecting point $x'_b$, which may be the same as $x_b$, or not (see Figure~\ref{fig:layers_gamma}). This way, we have a 1-to-1 relation between $\gamma^-$ and the intersecting points. For even higher order accuracy, layers of $\gamma^-$ can still be defined recursively and intersecting points can be selected accordingly. 
Then we recall the compatibility condition and impose it at the intersection points that correspond to the second layer:
\begin{align}\label{eqn:bc2_O4}
    \sum_{x_{j,k}\in \zeta} u_{j,k}\Delta\phi^{(3,3)}_{j,k}(x'_b) = f(x'_b).
\end{align}
Now, \eqref{eqn:bc1_O4} together with \eqref{eqn:bc2_O4} will account for $|\gamma^-|$ number of points. 
Normally the outer layer $\gamma^-_2$ will include more points than the inner layer $\gamma^{-}_1$, which implies that we will have different number of equations between \eqref{eqn:bc1_O4} and \eqref{eqn:bc2_O4} and they might be imposed at different points, unlike the 1D case where the number of equations are equal between Dirichlet boundary condition and compatibility conditions, which are imposed at the same boundary point.

\begin{remark}
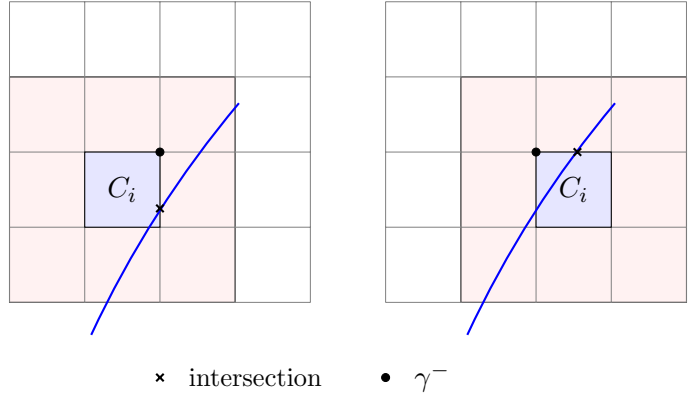
\begin{figure}[H]
    \centering
    \begin{tikzpicture}[scale=1]
    \filldraw[fill=red!5] (-2,-1) rectangle (1,2);
    \draw[step=1cm,gray,very thin] (-2,-1) grid (2,3);
    \filldraw[fill=blue!10] (-1,0) rectangle (0,1);
    \draw[blue,thick] (1.05,1.65) arc (140:155:14);
    \filldraw (0,1) circle (1.5pt);
    \draw (-0.5,0.5) node {$C_i$};
    \draw[mark size=+2pt,thick] plot[mark=x] coordinates {(0,0.25)};
    \filldraw[fill=red!5] (2+2,-1) rectangle (5+2,2);
    \draw[step=1cm,gray,very thin] (2+1,-1) grid (5+2,2+1);
    \filldraw[fill=blue!10] (3+2,0) rectangle (4+2,1);
    \draw[blue,thick] (1.05+3+2,1.65) arc (140:155:14);
    \filldraw (3+2,1) circle (1.5pt);
    \draw (3.5+2,0.5) node {$C_i$};
    \draw[mark size=+2pt,thick] plot[mark=x] coordinates {(3.55+2,1)};
    \draw[mark size=+2pt,thick] plot[mark=x] coordinates {(0,-2)};
    \draw (0.25,-2) node[right] {\small intersection};
    \filldraw (3,-2) circle (1.5pt);
    \draw (3.25,-2) node[right] {$\gamma^-$};
    \end{tikzpicture}
    \caption{Intersection points and their support cells}
    \label{fig:support_cell_O4}
\end{figure}
Depending on which intersection point we choose for the $\gamma^-$ point, the supporting cells will be different, as can be see in Figure~\ref{fig:support_cell_O4}. The point sets $\eta$ and $\gamma$ will be determined by the supporting cells as in the O2 case.
\end{remark}


Lastly, we are ready to present the set of boundary equations from which we will solve for the unknown density $u_{\gamma}$:
\begin{align}\label{eqn:boundary_eq_O4}
    \left(
    \begin{array}{c c c}
        0 & \multicolumn{2}{c}{I-P_{\zeta^+}} \\
        \multicolumn{3}{c}{\Phi^{(3,3)}(x_b)} \\
        \multicolumn{3}{c}{\Delta\Phi^{(3,3)}(x'_b)} \\
        \multicolumn{3}{c}{R}\\
    \end{array}
    \right)
    \begin{pmatrix}
        u_\eta\\
        u_\gamma\\
        u_\omega
    \end{pmatrix}
=
    \begin{pmatrix}
        G_hf_{\zeta^+}\\
        g(x_b)\\
        f(x'_b)\\
        0
    \end{pmatrix},
\end{align}
where $\Phi^{(3,3)}$ and $\Delta\Phi^{(3,3)}$ are matrices assembled from the basis functions and its laplacian evaluated at intersecting points $x_b$ and $x'_b$, respectively. Again, \eqref{eqn:boundary_eq_O4} gives a square system, where we can obtain the density $u_\gamma$ and apply it to the discrete generalized Green's formula \eqref{eqn:green} to get the approximated solution.

The developed numerical schemes are summarized into the following Algorithm \eqref{alg:local_dpm}.
\begin{algorithm}
\algFontSize 
\caption{Difference Potentials Method based on local basis functions}
\begin{algorithmic}[1]
    \State Embed $\Omega$ into an auxiliary domain $D$, and discretize $D$ with uniform Cartesian meshes from which $N^0,N^\pm,M^0,M^\pm,\gamma,\gamma^{\pm}$ are defined, given $n$-th order accurate centered finite difference stencils.
    \State Obtain intersection points $x_b$ that corresponds to $\gamma^-$.
    \State Construct supporting points $\eta$ and $\omega$ based on the intersection points and assemble $\zeta:=\eta\cup\gamma\cup\omega$.
    \For {$j=1,\dots,|\zeta|$}
        \State Obtain the $j$-th column of $P_\zeta$ by computing $G_{nh}X_{M^-}L_{nh}e_j$ where $e_j$ is an unit density at $x_j$.
    \EndFor
    \State Restrict $P_\zeta$ to $P_{\zeta^+}$.
    \State Evaluate basis functions (and their derivatives if necessary) at $x_b$ for all points in $\zeta$.
    \State Assemble the boundary equations into linear systems \eqref{eqn:boundary_eq_O2} or \eqref{eqn:boundary_eq_O4} and solve for density $u_\gamma$.
    \State Use the generalized discrete Green's formula \eqref{eqn:green} to obtain the approximation of the solutions.
\end{algorithmic}
\label{alg:local_dpm}
\end{algorithm}

\subsection{Fast Poisson solvers}
The major computational cost in the developed local basis approach lies in the construction of difference potentials. We need to solve the following system
\begin{align}\label{eqn:aux_unit}
L_{nh} u_h = X_{M^-}L_{nh}e_{j},\quad j=1,2,\dots,|\gamma|
\end{align}
subject to some auxiliary boundary conditions $|\gamma|$ times, where $e_{j}$ denotes a unit density with unit value at $x_{j}$ and 0 value at all other grid points. Here, $j$ denotes the index of a point in the set $\gamma$. Ideally, we would like to solve the auxiliary problem~\eqref{eqn:aux_unit} only once and translate the solutions to proper positions. However, such property only applies when $\sigma$ is sufficiently large such that the corresponding Green's function is localized.

Below we introduce a fast Fourier transform (FFT) based solver that applies to any order of spatial discretizations with centered finite difference stencils, in contrast with fast solvers based on compact finite difference stencils, which is limited to fourth order accuracy.

\paragraph{Second order accuracy}
The auxiliary problem defined on the rectangular domains can be solved using fast Poisson solvers based on FFT if it is of constant coefficients. In the second order case, homogeneous Dirichlet boundary condition leads to the following tridiagonal or block tridiagonal coefficient matrix:
\begin{subequations}
\begin{align}
K_{1d,2h} &= \mbox{diag}\left(1,-2,1\right),\\
K_{2d,2h} &= K_{1d,2h}\otimes I + I\otimes K_{1d,2h},\\
K_{3d,2h} &= K_{1d,2h}\otimes I \otimes I + I\otimes K_{1d,2h} \otimes I + I\otimes I\otimes K_{1d,2h}.
\end{align} 
\end{subequations}
It is well known that the coefficient matrix $K_{1d,2h}$ of size $N\times N$ can be diagonalized
\begin{align}
K_{1d,2h} = S^{-1}\Lambda^{(2h)} S
\end{align}
where $\Lambda^{(2h)}$ is a diagonal matrix with $j$-th entry given by
\begin{align}
\Lambda^{(2h)}_{jj}=2\left(\cos\left(\frac{j\pi}{N+1}\right)-1\right).
\end{align}
$S$ is a symmetric sine matrix with the $(j,k)$-th entry given by
\begin{align}
S_{jk} = \sin\left(\frac{jk\pi}{N+1}\right).
\end{align}
and its inverse is readily given by
\begin{align}
S^{-1} = \frac{2}{N+1}S
\end{align}

Diagonalization of coefficient matrices in higher dimensions can be obtained using the Kronecker product. In 2D, we have
\begin{subequations}
\begin{align}
K_{2d,2h} &= K_{1d,2h}\otimes I + I\otimes K_{1d,2h}\\
& = (S^{-1}\Lambda^{(2h)} S)\otimes (S^{-1}IS) + (S^{-1}IS)\otimes (S^{-1}\Lambda^{(2h)} S)\\
& = (S^{-1}\otimes S^{-1})(\Lambda^{(2h)}\otimes I+I\otimes \Lambda^{(2h)})(S\otimes S) \\
& = \left(\frac{2}{N+1}\right)^2(S\otimes S)(\Lambda^{(2h)}\otimes I+I\otimes \Lambda^{(2h)})(S\otimes S) 
\end{align}
\end{subequations}
and in 3D
\begin{subequations}
\begin{align}
K_{3d,2h} &= K_{1d,2h}\otimes I \otimes I + I\otimes K_{1d,2h} \otimes I + I\otimes I\otimes K_{1d,2h}\\
&=\left(\frac{2}{N+1}\right)^3(S\otimes S \otimes S)(\Lambda^{(2h)}\otimes I\otimes S+I\otimes \Lambda^{(2h)} \otimes I + I\otimes I \otimes \Lambda^{(2h)})(S\otimes S \otimes S)
\end{align}
\end{subequations}

To solve the discrete Poisson equation $K_{1d,2h}u_h = f_h$, we can apply the diagonalization. For example, in 1D, we have
\begin{align}
S^{-1}\Lambda S u_h = f_h \Rightarrow u_h = S^{-1}\Lambda^{-1}Sf = \frac{2}{N+1}S(\Lambda^{(2h)})^{-1}Sf_h
\end{align}
Higher dimensions follow similarly.

The exact steps of using FFT to solve the auxiliary problems on square domains are summarized in Algorithm~\ref{alg:fst} below.
\begin{algorithm}
\algFontSize 
\caption{FFT based fast Poisson solver for $K_{nd,2h}u_h = f_h$ in Dimension $n$}
\begin{algorithmic}[1]
    \If{$n=1$} \Comment Fast Sine Transform $f_h$ in $n$ dimensions
        \State $x = Sf_h$
    \ElsIf{$n=2$}
        \State $x = S\otimes S f_h$
    \ElsIf{$n=3$}
        \State $x = S\otimes S\otimes S f_h$
    \EndIf
    \If{$n=1$} \Comment Solve in Fourier Space
        \State Pointwise space loop: $y_j = 2/(N+1)x_j/\Lambda^{(2h)}_{jj}$
    \ElsIf{$n=2$}
        \State Pointwise space loop: $y_{jk} = [2/(N+1)]^2x_{jk}/(\lambda^{(2h)}_{jj}+\Lambda^{(2h)}_{kk})$
    \ElsIf{$n=3$}
        \State Pointwise space loop: $y_{jkl} = [2/(N+1)]^3x_{jkl}/(\Lambda^{(2h)}_{jj}+\Lambda^{(2h)}_{kk}+\Lambda^{(2h)}_{ll})$
    \EndIf
    \If{$n=1$} \Comment another Fast Sine Transform 
        \State $u_h = Sy$
    \ElsIf{$n=2$}
        \State $u_h = S\otimes S y$
    \ElsIf{$n=3$}
        \State $u_h = S\otimes S\otimes S y$
    \EndIf
\end{algorithmic}
\label{alg:fst}
\end{algorithm}

\paragraph{High order accuracy}

A key observation in \cite{FENG2020109391} for high-order finite difference discretizations was that antisymmetric boundary condition ($u(a-x)=-u(a+x)$ for any $x$ where $a$ is a boundary point) together with homogeneous Dirichelet boundary condition will allow a sine series expansion of the underlying grid functions. In the difference potentials framework, the choice of auxiliary boundary condition is arbitrary, hence we will use the antisymmetric boundary conditions to allow fast Poisson solvers based on FFT for high order accurate finite difference discretizations. Below we will illustrate this idea using the fourth order case, and higher order cases follow closely.

With the anti-symmetric boundary conditions and homogeneous Dirichlet BC, fourth order finite difference discretizations in 1D gives the following coefficient matrix
\begin{align}
K_{1d,4h} = \mbox{diag}\left(-\frac{1}{12},\frac{4}{3},-\frac{5}{2},\frac{4}{3},-\frac{1}{12}\right) + \frac{1}{12}E_1 + \frac{1}{12}E_N
\end{align}
where the entries of $E_1$ and $E_N$ assume zero value everywhere except $E_1(1,1)=1$ and $E_N(N,N)=1$. It can be shown that $K_{1d,4h}$ shares the same eigenvectors as $K_{1d,2h}$, i.e.,
\begin{align}
K_{1d,4h}S = S\Lambda^{(4h)}
\end{align}
but with a different eigenmatrix, which is given by
\begin{align}
\Lambda^{(4h)}_{jj} = -\frac{(\cos(j\pi/(N+1))-1)(\cos(j\pi/(N+1))-7)}{3}
\end{align}
This means Algorithm~\ref{alg:fst} can be employed for fourth order (and higher order) accurate finite difference discretizations as well, if one replaces $\Lambda^{(2h)}$ with $\Lambda^{(4h)}$. Sixth order or eighth order accurate versions can be designed similarly, with difference only in the eigenvalues, which can also be found in \cite{FENG2020109391}.

\section{Justification of difference potentials in max norm}\label{sec:convergence}
In this section, we follow \cite{ryaben2001method} and justify the use of difference potentials to approximate the singular boundary integrals in the solution of elliptic equations in $\infty$-norm, where H\"{o}lder's norm was used in \cite{ryaben2001method}. Some of the well-known results are given here without proof. For simplicity, we consider Poisson's equation in an arbitrary domain $\Omega$:
\begin{align}\label{eqn:poisson}
    L u = f, \quad x\in \Omega,
\end{align}
subject to the boundary condition
\begin{align}
    \ell u = g,\quad x\in \Gamma. \label{eqn:bc_in}
\end{align}
Assume the Cauchy data of the solution $\bm{u}_\Gamma := (u,\frac{\partial u}{\partial n})|_\Gamma$ is known at the boundary, then the solution can be expressed from the Green's formula
\begin{align}
    u = \int_\Gamma u\frac{\partial G}{\partial n}-\frac{\partial u}{\partial n}Gds_y + \int_\Omega Gf dy
\end{align}
or in the operator form
\begin{align}
    u = P_{\Omega\Gamma}\bm{u}_\Gamma + Gf
\end{align}
where $P_{\Omega\Gamma}\bm{u}_\Gamma$ is the potential defined by 
\begin{align}\label{eqn:potential1}
    P_{\Omega\Gamma}\bm{u}_\Gamma = \int_\Gamma u\frac{\partial G}{\partial n}-\frac{\partial u}{\partial n}Gds_y = u - Gf= u - GL u
\end{align}

\begin{remark}
It can be seen from the definition of the differential potentials $P_{\Omega\Gamma}\bm{u}_\Gamma$ that it only depends on the Cauchy data $\bm{u}|_\Gamma$, i.e., $u_1 - GL u_1 = u_2 - GL u_2$ if $\bm{u}_1|_\Gamma=\bm{u}_2|_\Gamma$, which will provide us the convenience to choose the function to compute the potentials.
\end{remark}

Due to the lack of knowledge of the Green's function $G$ for a general domain $\Omega$, the differential potentials in \eqref{eqn:potential1} is not obtainable. Instead, we introduce another equivalent definition in an auxiliary domain. To do this, we first introduce an auxiliary problem defined in a larger computationally simple domain $D\supset\Omega$,
\begin{align}\label{eqn:aux}
    L \tilde{u} = \tilde{f}, \quad x\in D
\end{align}
subject to some boundary condition for $(x,y)\in\partial D$, where $\tilde{u}$ and $\tilde{f}$ are some smooth extension of $u$ and $f$ from $\Omega$ to $D$, such that the Cauchy data of $\tilde{u}(x,y)$ coincides with the Cauchy data of $u(x,y)$.

\begin{definition}[Differential Potentials]
The differential potentials can be computed as
\begin{align}\label{eqn:potential2}
    P_{D\Gamma}\bm{u}_\Gamma = \tilde{u} - \tilde{G}[X_{D}L \tilde{u}]
\end{align}
given the vector density $\bm{u}_\Gamma$.
\end{definition}
Clearly, we have
\begin{align}
    P_{\Omega\Gamma}\bm{u}_\Gamma  = Tr_\Omega P_{D\Gamma}\bm{u}_\Gamma
\end{align}
since they both satisfy the Laplace equation and the same boundary condition at $\partial\Omega$.

For the discrete version, we first discretize the auxiliary domain $D$ using uniform Cartesian meshes, and introduce point sets $M^\pm$, $M^0$, $N^\pm$, $N^0$ and $\gamma$ using $p$-th order accurate centered finite difference schemes as defined in Section~\ref{sec:scheme}. Then Possion's equation \eqref{eqn:poisson} can be discretized as
\begin{align}
    L_h u_{jk} = f_{jk},\quad x_{jk} \in M^+
\end{align}
subject to boundary condition $u = u_\gamma$ for point set $\gamma$. Here we assume $u_\gamma$ is known from the boundary condition. The difference potentials for a given density $u_\gamma$ is computed as
\begin{align}\label{eqn:diff_potential1}
    P_{N^+\gamma}u_\gamma = u_{jk} - G_h[X_{M^+}L_h u_{jk}]
\end{align}
where where the grid function satisfies the trace condition $Tr_\gamma u_{jk} = u_\gamma$. Again, difference potentials are not actually computable by \eqref{eqn:diff_potential1}, thus we first discretize he auxiliary problem \eqref{eqn:aux} as
\begin{align}
    L_h \tilde{u}_{jk} = \tilde{f}_{jk},\quad x_{jk} \in M^0
\end{align}
subject to homogeneous boundary conditions plus antisymmetric boundary condition if necessary, which we will use to do the actual computation of the difference potentials. Here we will not differentiate between $L_h$ on $M+$ or $M^0$, since there will be no risk of confusion.

\begin{definition}[Difference Potentials]
The discrete analog of the differential potentials \eqref{eqn:potential2} are difference potentials, which are computed as
\begin{align}\label{eqn:diff_potential2}
    P_{N^0\gamma}u_\gamma = \tilde{u}_{jk} - \tilde{G}_h[\chi_{M^+}\tilde{L}_h \tilde{u}_{jk}]
\end{align}
for any given density $u_\gamma$ and $Tr_\gamma\tilde{u}_{jk}=u_\gamma$.
\end{definition}
Using similar argument as in the continuous case, we may conclude that the difference potentials obtained from the auxiliary problem is the same as the one obtained in \eqref{eqn:diff_potential1}:
\begin{align}
    P_{N^+\gamma}u_\gamma = Tr_{N^+}P_{N^0\gamma}u_\gamma
\end{align}

Up till now, what we introduced in this section is standard in the difference potentials framework as in \cite{ryaben2001method}. With the differential and difference potentials operators defined, we are ready to introduce the following lemma.

\begin{lemma}\label{lemma:lap}
    Assume we have a function $v(x)$ that satisfies the auxiliary problem \eqref{eqn:aux} with homogeneous boundary condition (plus anti-symmetric boundary condition if necessary), and the trace condition $\bm{v}|_\Gamma=\bm{u}|_{\Gamma}$. Let the grid function be $v_{N^0}=Tr_{N^0}v(x)$, and the density be $v_\gamma = Tr_\gamma v_{N^0}$, then
    \begin{align}
        \left|\left|P_{N^0\gamma}v_\gamma-Tr_{N^0}P_{D\Gamma}\bm{v}_\Gamma\right|\right|_\infty\leq C h^{p}
    \end{align}
    where $p$ is the order of accuracy of $L_h$ for approximating $L$.
\end{lemma}

\begin{proof}
    First note that
    \begin{align}
        P_{N^0\gamma}v_\gamma = v_{N^0}-\tilde{G}_hX_{M^+}L_hv_{N^0}
    \end{align}
    and
    \begin{align}
        Tr_{N^0}P_{D\Gamma}\bm{v}_\Gamma = Tr_{N^0}\left[v_{D}-\tilde{G}X_{D}Lv\right]
    \end{align}
    and the difference
    \begin{subequations}
    \begin{align}
        w_{N^0} &= Tr_{N^0}P_{D\Gamma}\bm{v}_\Gamma-P_{N^0\gamma}v_\gamma\\
        &= \tilde{G}_hX_{M^+}L_hv_{N^0}-Tr_{N^0}\tilde{G}\chi_{D}Lv
    \end{align}
    \end{subequations}
    The first term is a solution to
    \begin{align}
        L_h z_{N^0} = X_{M^+}L_hv_{N^0}, \quad x_{jk}\in M^0
    \end{align}
    subject to homogeneous boundary condition
    and the second term is the trace of the solution to
    \begin{align}
        L z = X_{D}Lv, \quad x\in D
    \end{align}
    subject to same boundary condition.

    By assumption on the accuracy of $L_h$ and the discrete PDE, we have
    \begin{align}
        L_hz_{N^0} = X_{M^+}[L_hv_{N^0}]= X_{M^+}[Tr_{M^0}Lv + h^p\phi_{M^0}],\quad x_{jk}\in M^0
    \end{align}
    for some grid function $phi_{M^0}$ that satisfies $||\phi_{M^0}||_\infty\leq C_1$.
    For the second term, we have
    \begin{align}
        L_h[Tr_{N^0}z]= Tr_{M^0}Lz + h^p\psi_{M^0}= Tr_{M^0}[X_{D}Lv] + h^p\psi_{M^0} ,\quad x_{jk}\in M^0
    \end{align}
    for some grid function $phi_{M^0}$ that satisfies $||\psi_{M^0}||_\infty\leq C_2$.
    Then the difference $w_{N^0}$ satisfies
    \begin{subequations}
    \begin{align}
        L_h w_{N^0} &= L_hz_{N^0}-L_h[Tr_{N^0}z]\\
        &=X_{M^+}[Tr_{M^0}Lv + h^p\phi_{M^0}]-Tr_{M^0}[X_{D}Lv] - h^p\psi_{M^0}\\
        &=X_{M^+} h^p\phi_{M^0}-h^p\psi_{M^0}\\
        &=h^p(X_{M^+}\phi_{M^0}-\psi_{M^0})
    \end{align}
    \end{subequations}
    Lastly, the difference can be bounded as
    \begin{align}
        ||w_{N^0}||_\infty\leq C||L^{-1}_h||_{\infty}h^p
    \end{align}
    which concludes our proof, since $L^{-1}_h$ is uniformly bounded in $\infty$ norm (see, e.g. \cite{hackbusch1981regularity,hackbusch2017elliptic}).
\end{proof}

\begin{lemma}\label{lemma:ugamma}
Suppose that $u$ satisfies $Lu=f$ in a domain $D$ with homogeneous boundary condition, and $u_h$ satisfies $L_hu_h=f_h$ on $D_h$ where the subscript denotes their discrete counterparts. Assume $L_h$ approximates $L$ with accuracy $h^p$, and $u_h$ approximates $u$ with accuracy $h^q$. Let the difference of traces of $u$ and $u_h$ on $\gamma$ be
\begin{align}
    v_\gamma = Tr_{\gamma}u-Tr_{\gamma}u_h,
\end{align}
then the associated difference potentials obtained with the density $v_\gamma$ satisfies the following estimate:
\begin{align}
    ||P_{N^+\gamma}v_\gamma||_{\infty}\leq Ch^{\min(p,q)}.
\end{align}
\end{lemma}

\begin{proof}
    We first extend $v_\gamma$ to $N^0$ by 0 to get $v_{N^0}$, then the discretization of $v_{N^0}$ satsifies 
    \begin{align}
        ||L_hv_{N^0}||_\infty\leq ||L_h(Tr_{N^0}u-u_h)||_\infty\leq C(h^p+||Tr_{N^0}u-u_h||_{\infty}) \leq Ch^{\min(p,q)}.
    \end{align}
    following the proof of the interior estimate in \cite[Theorem 4.1]{Thom_e_1968}.
    
    Note that the difference potentials $P_{N^+\gamma}v_\gamma$ are computed as
    \begin{align}
        V_{N^0} = v_{N^0}-\tilde{G}_h\chi_{M^+}L_hv_{N^0},\quad P_{N^+\gamma}v_\gamma = Tr_{N^+}V_{N^0}
    \end{align}
    which means $V_{N^0}$ satisfies the following system,
    \begin{align}
        L_h V_{N^0} =  L_hv_{N^0}- L_h\tilde{G}_h\chi_{M^+}L_hv_{N^0}=\chi_{M^-}L_hv_{N^0},
    \end{align}
    subject to homogeneous Dirichlet boundary condition.

    Thus, we have
    \begin{align}
        ||V_{N^0}||_{\infty} = ||L_h^{-1}\chi_{M^-}L_hv_{N^0}||_{\infty}\leq ||L_h^{-1}||_{\infty}||\chi_{M^-}||_{\infty}||L_hv_{N^0}||_{\infty}\leq Ch^{\min(p,q)}.
    \end{align}
    which implies the difference potentials $P_{N^+\gamma}$ also satisfy the estimate.
\end{proof}


With the two lemmas \ref{lemma:lap} and \ref{lemma:ugamma}, we are ready to give the theorem that establishes the approximation of differential potentials (the singular integrals in the Poisson's equation case)

\begin{theorem}\label{thm:conv}
    If the density $u_\gamma$ approximates the the trace of $u$ in gamma (i.e. $Tr_{\gamma}u$) with accuracy of order $q$,
    \begin{align}
        ||u_\gamma-Tr_{\gamma}u||_{\infty} \leq Ch^q
    \end{align}
    then the difference potentials approximate the differential potentials with order of accuracy $\min(p,q)$, i.e.,
    \begin{align}
        \left|\left|P_{N^+\gamma}u_\gamma - Tr_{N^+}P_{D\Gamma}\bm{u}_{\Gamma}\right|\right|_{\infty} \leq C h^{\min(p,q)}
    \end{align}
    provided the discrete operator $L_h$ approximates the continuous $L$ with order of accuracy $p$.
\end{theorem}

\begin{proof}
    By triangle inequality,
    \begin{subequations}
    \begin{align}
        &\left|\left|P_{N^+\gamma}u_\gamma - Tr_{N^+}P_{D\Gamma}\bm{u}_{\Gamma}\right|\right|_{\infty} \\
        \leq& \left|\left|P_{N^+\gamma}Tr_{\gamma}u - Tr_{N^+}P_{D\Gamma}\bm{u}_{\Gamma}\right|\right|_{\infty} + \left|\left|P_{N^+\gamma}u_\gamma - P_{N^+\gamma}Tr_{\gamma}u\right|\right|_{\infty}  \\
        \leq&  Ch^{\min(p,q)}
    \end{align}
    \end{subequations}
    where the last two inequalities are by Lemmas \ref{lemma:lap} and \ref{lemma:ugamma}, respectively.
\end{proof}

\begin{remark}
Theorem~\ref{thm:conv} implies that the accuracy of the difference potentials is limited by both the accuracy of the spatial discretizations and the accuracy of the density. In practice with the local basis functions, the accuracy of the density $u_\gamma$ will not exceed that of the spatial discretizations, thus the limit factor will be $q$. For Dirichlet BC, we can have $q=p$, but with Neumann BC, $q=p-1$, which is also observed in the numerical section.

\end{remark}

Now the only remaining question is whether the density $u_\gamma$ we obtained using the boundary equations introduced in Section~\ref{sec:scheme} is accurate enough. Below we show an example of Dirichlet BC in 1D that he density $u_\gamma$ obtained from solving boundary equations \eqref{eqn:boundary-eqs} is of second order accuracy.

\begin{theorem}\label{thm:O2_density_accuracy}
The density $u_\gamma$ obtained from solving boundary equations \eqref{eqn:boundary-eqs} gives $\mathcal{O}(h^2)$ approximation to $Tr_{\gamma}u$ in $\infty$ norm, i.e.,
\begin{align}
||u_\gamma-Tr_\gamma u||_{\infty} \leq C h^2
\end{align}
where $u$ is a solution to the continuous auxiliary problem.
\end{theorem}

\begin{proof}
We will use the Poisson's equation as an example. 
The boundary equation 
\begin{align}
u_{\gamma^+}-P_{\gamma^+}u_{\gamma} &= G_hf_{\gamma^+}
\end{align}
is equivalent to the discrete PDE defined on $M^+$, i.e.,
\begin{align}
\frac{u_{j-1}-2u_j+u_{j+1}}{h^2} = f_j,\quad x_j\in M^+.
\end{align}
Without loss of generality, we will assume $j$ goes from 1 to $N$ with $j=1,N$ outside of $M^+$ and $j=2,\dots,N-1$ inside of $M^+$.
\begin{figure}[H]
\centering
\begin{tikzpicture}[scale=1.5]
\draw[thick] (-0.5,0) -- (4.5,0);
\filldraw (0.0,0) circle (1pt);
\draw (0,-0.1) node[below] {\footnotesize 2};
\filldraw (0.5,0) circle (1pt);
\filldraw (1.0,0) circle (1pt);
\filldraw (1.5,0) circle (1pt);
\filldraw (2.0,0) circle (1pt);
\filldraw (2.5,0) circle (1pt);
\filldraw (3.0,0) circle (1pt);
\filldraw (3.5,0) circle (1pt);
\filldraw (4.0,0) circle (1pt);
\draw (4.0,-0.1) node[below] {\footnotesize $N-1$};

\draw (-0.5,0) circle (1pt);
\draw (-0.5,-0.1) node[below] {\footnotesize 1};

\draw (4.5,0) circle (1pt);
\draw (4.5,-0.1) node[below] {\footnotesize $N$};

\draw[mark size=+1.5pt,thick] plot[mark=x] coordinates {(-0.25,0)};
\draw[mark size=+1.5pt,thick] plot[mark=x] coordinates {(4.25,0)};
\draw (-0.25,0.1) node[above] {$a$};
\draw (4.25,0.1) node[above] {$b$};

\end{tikzpicture}
\caption{An illustration of point set $M^+$ (solid), $\gamma^-$ (circle) and boundary points $a,b$}
\label{fig:proof_full_grids}
\end{figure}
We will use points in Figure~\ref{fig:proof_full_grids} and assume $a$ and $b$ can lie anywhere between $[x_1,x_2]$ and $[x_{N-1},x_{N}]$, respectively. Now the boundary equations will be equivalent to the following system
\begin{subequations}
\begin{align}
\frac{u_{j-1}-2u_j+u_{j+1}}{h^2} &= f_j,\quad j=2,\dots,N-1\\
u_1\phi_1(a) + u_2\phi_2(a) &= u_a\\
u_{N-1}\phi_{N-1}(b) + u_{N}\phi_{N}(b) &= u_b
\end{align}
\end{subequations}
where the coefficient matrix can be written as
\begin{align}
M = \begin{pmatrix}
\phi_1(a) & \phi_2(a) & & &\\
1/h^2 & -2/h^2 & 1/h^2 & &\\
&\ddots &\ddots &\ddots & \\
& & 1/h^2 & -2/h^2 & 1/h^2\\
& & & \phi_{N-1}(b) & \phi_{N}(b)
\end{pmatrix}
\end{align}
where the basis functions $\phi_1(a),\phi_2(a),\phi_{N-1}(b),\phi_{N}(b)$ will lie in $[0,1]$.

Since the local truncation errors are of order $\mathcal{O}(h^2)$, we only need to show $||M^{-1}||_{\infty}$ is bounded so as to show $u_\gamma$ obtained in \eqref{eqn:boundary-eqs} is of order $\mathcal{O}(h^2)$ in $\infty$-norm. Here the matrix $\infty$ norm is induced from the vector $\infty$ norm. To this end, we introduce a sequence of matrices that is of the same size of $M$: $A_0=diag(1/h^2, -2/h^2, 1/h^2)$, $A_1$ equals $A_0$ except the first row is replaced by the first row of $M$, and $A_2=M$. First, we have $||A^{-1}_0||_{\infty}\leq C$ and we can regard $A_1$ as a rank-1 perturbation of $A_0$ ($A_1 = A_0+u_0v^T_0$), and $A_2$ a rank-1 perturbation of $A_1$ ($A_2 = A_1 + u_1v_1^T$). We will start with the inverse of $A_1$ using the Sherman-Morrison-Woodbury formula, and use the same argument for $A_2$. First, we have
\begin{align}
A^{-1}_1 = A_0^{-1} - \frac{A^{-1}_0u_0v^T_0A^{-1}_0}{1+v^T_0A^{-1}_0u_0}
\end{align}
where $u_0=(1,0,\dots,0)^T$, $v_0 = (\phi_1(a)+2/h^2,\phi_2(a)-1/h^2,\dots,0)^T$. Now note the denominator can be explicitly computed as
\begin{subequations}
\begin{align}
1+v^T_0A^{-1}_0u_0 &= 1+(\phi_1(a)+2/h^2)[A_0^{-1}]_{11}+(\phi_2(a)-1/h^2)[A^{-1}_0]_{21}\\
&=\left\{1+2/h^2[A_0^{-1}]_{11}-1/h^2[A^{-1}_0]_{21}\right\}+\phi_1(a)[A_0^{-1}]_{11}+\phi_2(a)[A^{-1}_0]_{21}\\
&=\phi_1(a)[A_0^{-1}]_{11}+\phi_2(a)[A^{-1}_0]_{21}
\end{align}
\end{subequations}
since $1=-2/h^2[A_0^{-1}]_{11}+1/h^2[A^{-1}_0]_{21}$ from the definition of inverse matrix $A_0^{-1}$. Then the numerator can also be explicitly expressed as
\begin{align}
A^{-1}_0u_0v^T_0A^{-1}_0 &= \begin{pmatrix}
[A^{-1}_0]_{11} (\phi_1(a)+2/h^2) & [A^{-1}_0]_{11}(\phi_2(a)-1/h^2) & 0 & \cdots & 0\\
[A^{-1}_0]_{21} (\phi_1(a)+2/h^2) & [A^{-1}_0]_{21}(\phi_2(a)-1/h^2) & 0 & \cdots & 0\\
\vdots & \vdots & \vdots & \ddots & \vdots\\
[A^{-1}_0]_{N1} (\phi_1(a)+2/h^2) & [A^{-1}_0]_{N1}(\phi_2(a)-1/h^2) & 0 & \cdots & 0\\
\end{pmatrix}A^{-1}_0 
\end{align}
Then each entry of $A^{-1}_0u_0v^T_0A^{-1}_0$ can be expressed as
\begin{subequations}
\begin{align}
\{A^{-1}_0u_0v^T_0A^{-1}_0\}_{jk}&= [A^{-1}_0]_{j1} (\phi_1(a)+2/h^2)[A^{-1}_0]_{1k} + [A^{-1}_0]_{j1}(\phi_2(a)-1/h^2)[A^{-1}_0]_{2k}\\
&= [A^{-1}_0]_{j1}\left\{ (\phi_1(a)+2/h^2)[A^{-1}_0]_{1k} + (\phi_2(a)-1/h^2)[A^{-1}_0]_{2k} \right\}\\
&= [A^{-1}_0]_{j1}\left\{ \phi_1(a)[A^{-1}_0]_{1k} + \phi_2(a)[A^{-1}_0]_{2k}+2/h^2[A^{-1}_0]_{1k}-1/h^2[A^{-1}_0]_{2k} \right\}\\
&= [A^{-1}_0]_{j1}\left\{ \phi_1(a)[A^{-1}_0]_{1k} + \phi_2(a)[A^{-1}_0]_{2k}-\delta_{1k} \right\}
\end{align}
\end{subequations}
and its $\infty$-norm can be estimated as
\begin{subequations}
\begin{align}
||A^{-1}_0u_0v^T_0A^{-1}_0||_{\infty} &= \max_{j}\sum_{k=1}^N\Big|[A^{-1}_0]_{j1}\Big|\Big| \phi_1(a)[A^{-1}_0]_{1k} + \phi_2(a)A^{-1}_0]_{2k}-\delta_{1k} \Big|\\
&\leq \max_{j}\left|[A^{-1}_0]_{j1}\right|\left\{ \phi_1(a)\sum_{k=1}^N\left|[A^{-1}_0]_{1k}\right| + \phi_2(a)\sum_{k=1}^N\left|[A^{-1}_0]_{2k}\right|+\delta_{1k} \right\}\\
&\leq \max_{j}\left|[A^{-1}_0]_{j1}\right|\left\{ \phi_1(a)||A_{0}^{-1}||_{\infty} + \phi_2(a)||A_{0}^{-1}||_{\infty}+1 \right\}\\
&\leq \max_{j}\left|[A^{-1}_0]_{j1}\right|\left\{ [\phi_1(a)+\phi_2(a)]||A_{0}^{-1}||_{\infty} +1 \right\}\\
&= \max_{j}\left|[A^{-1}_0]_{j1}\right|\left(||A_{0}^{-1}||_{\infty} +1 \right)
\end{align}
\end{subequations}
The last equality holds since $\phi_1(a)+\phi_2(a)=1$.
Now the inverse of $A_1$ can be estimated as
\begin{subequations}
\begin{align}
||A_1^{-1}||_{\infty} &\leq ||A^{-1}_0||_{\infty}+ \dfrac{\max_{j}\left|[A^{-1}_0]_{j1}\right|\left(||A_{0}^{-1}||_{\infty} +1 \right)}{\phi_1(a)\left|[A_0^{-1}]_{11}\right|+\phi_2(a)\left|[A^{-1}_0]_{21}\right|}\\
&\leq ||A^{-1}_0||_{\infty}+ \dfrac{\max_{j}\left|[A^{-1}_0]_{j1}\right|\left(||A_{0}^{-1}||_{\infty} +1 \right)}{[\phi_1(a)+\phi_2(a)]\min\left(\left|[A_0^{-1}]_{11}\right|,\left|[A^{-1}_0]_{21}\right|\right)}\\
&= ||A^{-1}_0||_{\infty}+ \dfrac{\max_{j}\left|[A^{-1}_0]_{j1}\right|}{\min\left(\left|[A_0^{-1}]_{11}\right|,\left|[A^{-1}_0]_{21}\right|\right)}\left(||A_{0}^{-1}||_{\infty} +1 \right)
\end{align}
\end{subequations}
which does not depend on $\phi_1(a)$ and $\phi_2(a)$, or the boundary point $a$.

We claim that $[A^{-1}_0]_{11}$ and $[A^{-1}_0]_{21}$ will not be 0. For example, if $[A^{-1}_0]_{11}=0$, this leads to $[A^{-1}_0]_{21}=[A^{-1}_0]_{31}=\dots=[A^{-1}_0]_{N1}=1$, which gives a contradiction in the last row. The case of $[A^{-1}_0]_{21}=0$ can be argued similarly. Now the denominator is finite and we need to check the numerator and we also claim that $\max_{j}\left|[A^{-1}_0]_{j1}\right|$ are bounded, as $||A^{-1}_0||_{\infty}$ is bounded. Thus it is safe to claim that $||A_1^{-1}||_{\infty}$ is also bounded by some constant $C$.

Next we can use the same argument for $A_2 = A_1 + u_1v_1^T$ with only difference in index in the computation, and we can conclude that $A_2$ or $M$ also satisfies
\begin{align}
||A_2^{-1}||_{\infty} = ||A^{-1}_1||_{\infty}+ \dfrac{\max_{j}\left|[A^{-1}_1]_{jN}\right|}{\min\left(\left|[A_1^{-1}]_{(N-1)N}\right|,\left|[A^{-1}_1]_{NN}\right|\right)}\left(||A_{1}^{-1}||_{\infty} +1 \right)\leq C
\end{align}
for some constant $C$, which concludes our proof.

\end{proof}

\begin{remark}
The proof outlined in Theorem~\ref{thm:O2_density_accuracy} depends on the non-negativity of the O2 basis functions $\phi_1$, $\phi_2$, $\phi_{N-1}$ and $\phi_N$, and the O4 basis functions fail to satisfy this condition (see Figure~\ref{fig:basis_and_deriv}), which means we need to find another approach to prove the accuracy of boundary equations \eqref{eqn:boundary-eqs-O4} and \eqref{eqn:boundary-eqs-compat}.
\end{remark}

\section{Numerical Results}\label{sec:numerics}
In this section, we validate the developed schemes numerically in 2D. Throughout this section, we will use the manufactured solution
\begin{align}
u(x,y) = \sin(x)\cos(y)
\end{align}
Forcing function $f(x,y)$ and boundary conditions will be obtained according to the manufactured solution. Unless specified otherwise, we will use the PDE $\Delta u -\sigma u = f$ with $\sigma=0$ and Dirichlet boundary conditions for testing the developed scheme. Other types of PDEs will not give significantly different results. The auxiliary boundary conditions are always chosen to be homogeneous, plus anti-symmetry for high-order schemes.

\subsection{An elliptical domain}
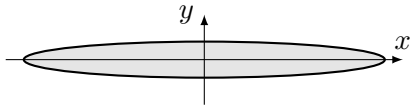
\begin{figure}[H]
\centering
\begin{tikzpicture}[scale=1.2]
  \draw[thick,fill=gray!40,fill opacity=.5] (0,0) ellipse (2cm and 0.2cm);
    \draw[->,>=latex] (-2.2,0)--(2.2,0) node[above] {$x$};
  \draw[->,>=latex] (0,-0.5)--(0,0.5) node[left] {$y$};
\end{tikzpicture}
\caption{An ellipse with $\alpha=10$}\label{fig:ellipse}
\end{figure}

The first example is an ellipse with aspect ratio $\alpha$.
\begin{align}
\psi(x,y) = x^2+\alpha^2y^2-1
\end{align}
which is embedded into an auxiliary domain $[-1.2,1.2]\times[-1.2,1.2]$. The extra length 0.2 is a random choice and we only need to make sure the auxiliary boundary is sufficient away from the ``real'' boundary to allow for a few grid widths for the auxiliary solution to decay to the homogeneous boundary condition. The ellipse with aspect ratio $\alpha=10$ is plotted in Figure~\ref{fig:ellipse}.

\begin{figure}[H]
    \centering
    \begin{subfigure}{0.45\textwidth}
        \centering
        \includegraphics[width=\textwidth, trim = 2cm 7cm 2cm 7.5cm]{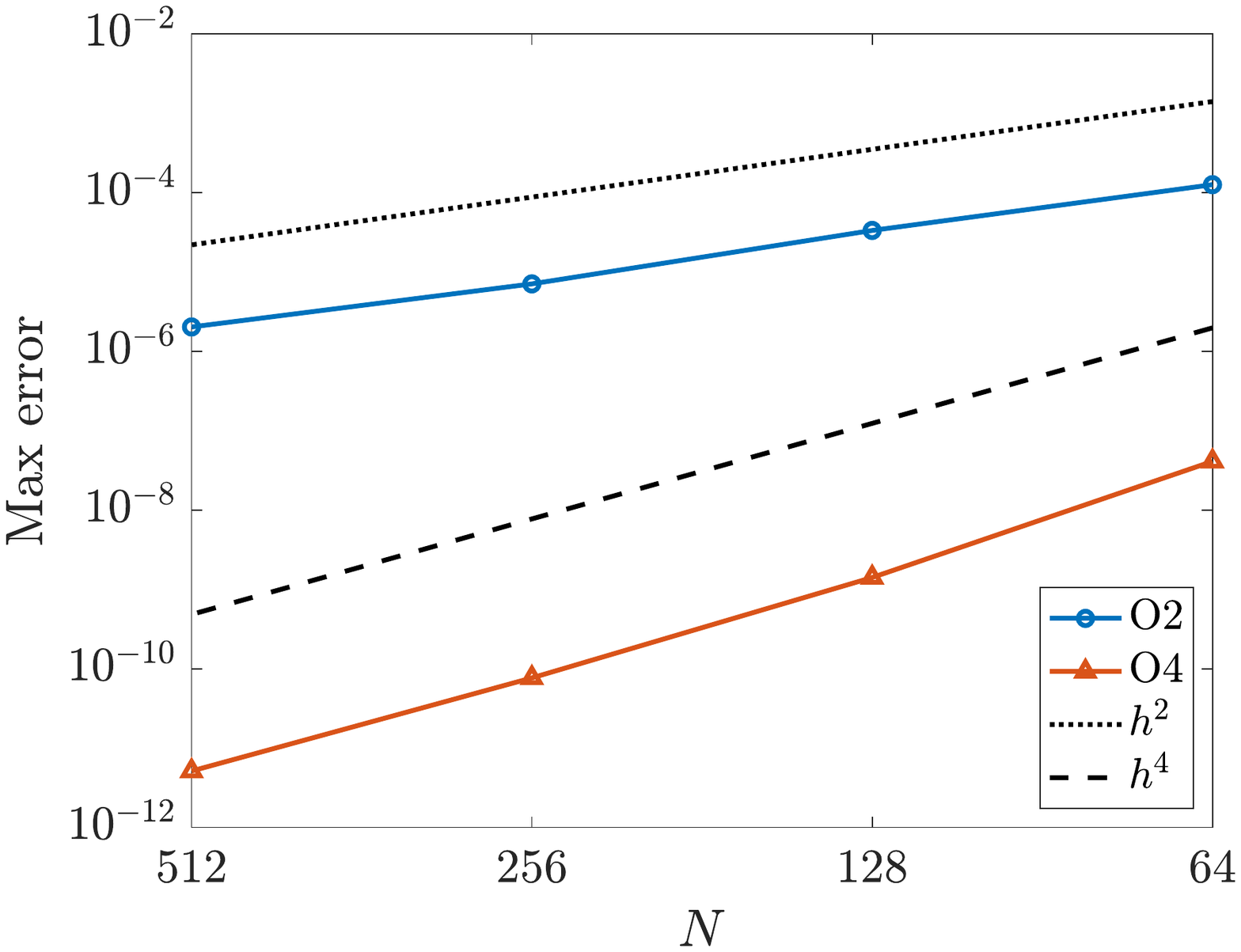}
        \caption{Convergence}
        \label{fig:ellipse_conv}
    \end{subfigure}
    \quad
    \begin{subfigure}{0.45\textwidth}
        \centering
        \includegraphics[width=\textwidth, trim = 2cm 7cm 2cm 7.5cm]{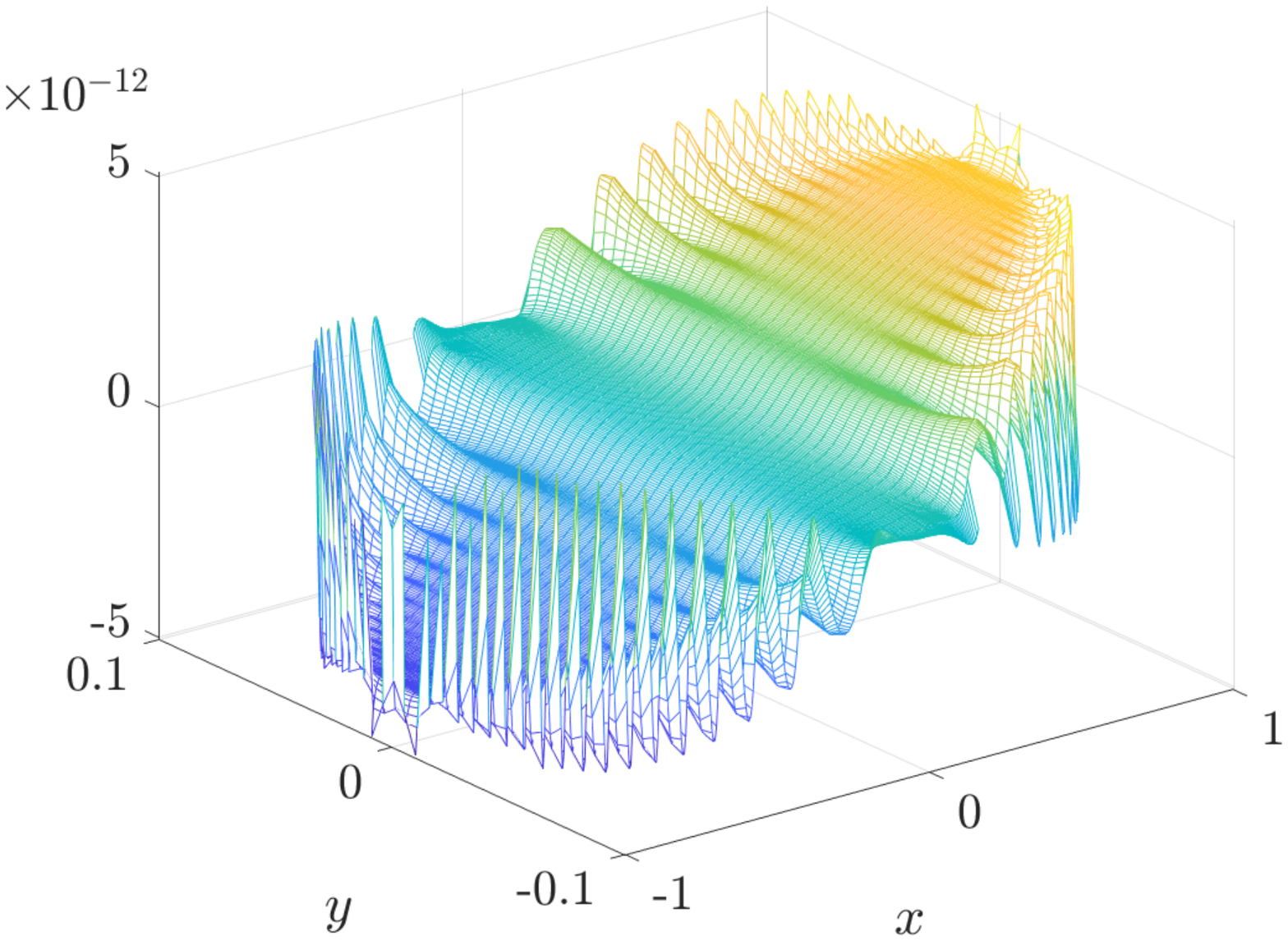}
        \caption{O4 errors $(N=512)$}
        \label{fig:ellipse_err}
    \end{subfigure}
    \caption{Convergence and errors for $\sigma=0$ and aspect ratio $\alpha=10$}\label{fig:ellipse_conv_err}
\end{figure}

The max errors and convergence are shown in Figure~\ref{fig:ellipse_conv_err}, where the second and fourth order convergence is observed. The errors exhibit random spikes near the boundary, which are dominating, and smooth oscillations of higher modes inside the domain. Those oscillations might benefit from some post-processing techniques such as those developed in \cite{docampo2020enhancing}.

\begin{figure}[H]
    \centering
    \begin{subfigure}{0.45\textwidth}
        \centering
        \includegraphics[width=\textwidth, trim = 2cm 7cm 2cm 7.5cm]{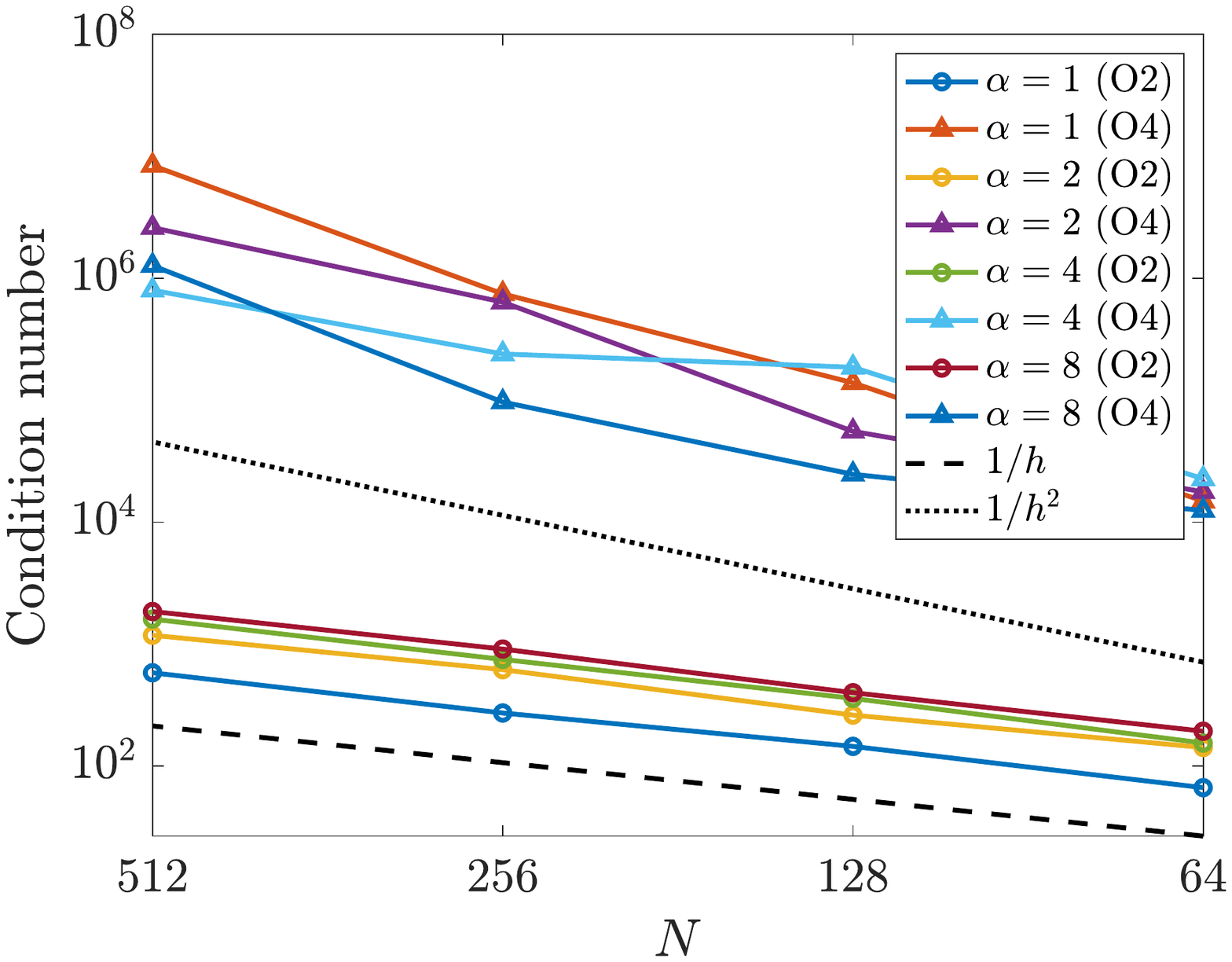}
        \caption{Different $\alpha$ for $\sigma=0$}
        \label{fig:ellipse_cond_diff_alpha}
    \end{subfigure}
    \quad
    \begin{subfigure}{0.45\textwidth}
        \centering
        \includegraphics[width=\textwidth, trim = 2cm 7cm 2cm 7.5cm]{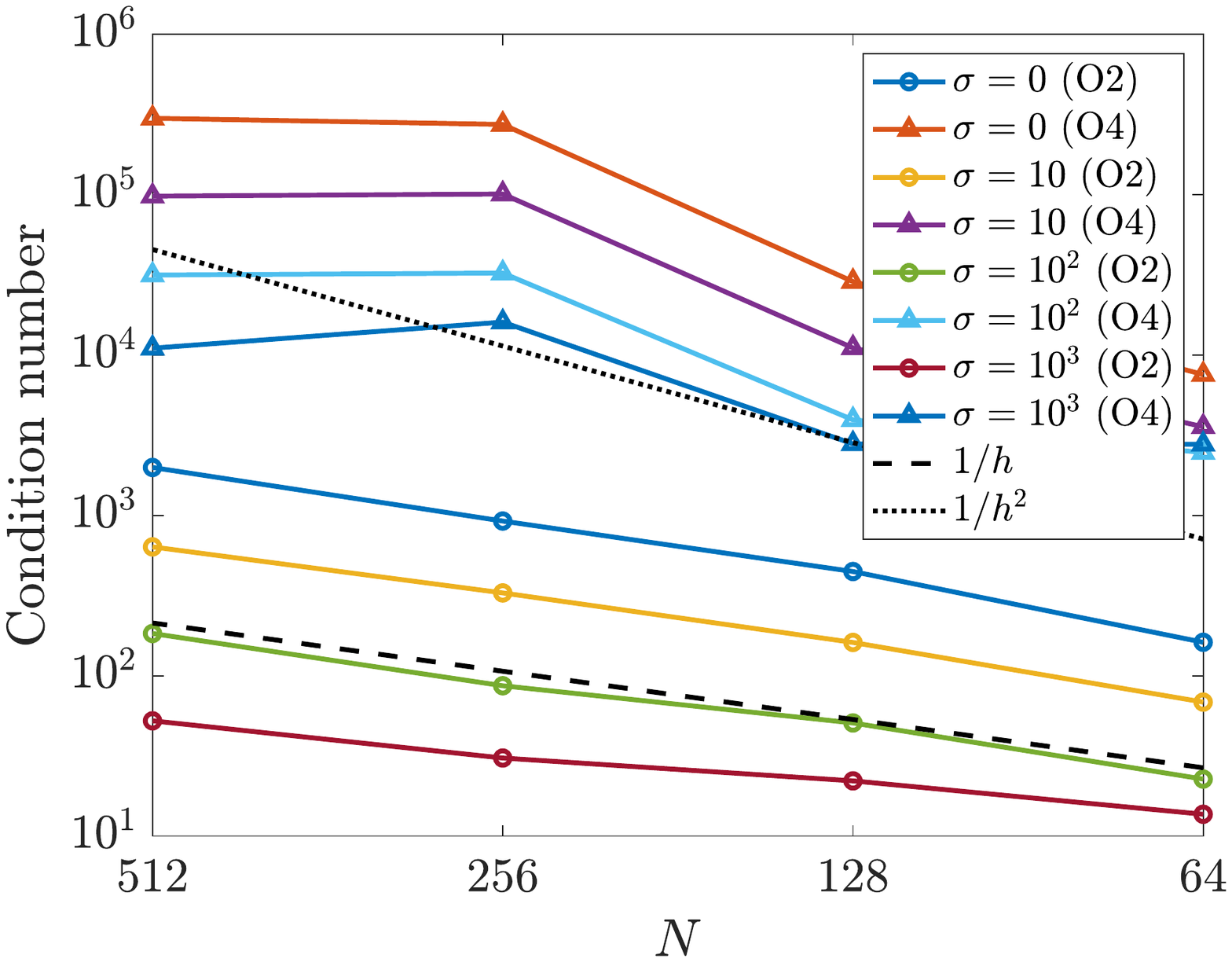}
        \caption{Different $\sigma$ for $\alpha=10$}
        \label{fig:ellipse_cond_diff_sigma}
    \end{subfigure}
    \caption{Condition numbers}\label{fig:ellipse_cond_diff}
\end{figure}

Next we look at different aspect ratios ($\alpha$ varies while $\sigma=0$) and also different PDEs ($\sigma$ varies while $\alpha=10$) in Figure~\ref{fig:ellipse_cond_diff}. It can be observed that for the same PDE in ellipses with different aspect ratios, the pattern of the growth of condition numbers is similarly. O2 schemes follows the linear growth line, while O4 exhibits an approximately quadratic growth rate. When we test different $\sigma$s using the same ellipse, we can see that the condition number is better with larger $\sigma$ values, which is desirable when this developed scheme is applied to time dependent problems when implicit time stepping is adopted.

\subsection{Robin boundary conditions}

\begin{figure}[H]
    \centering
    \begin{subfigure}{0.45\textwidth}
        \centering
        \includegraphics[width=\textwidth, trim = 2cm 7cm 2cm 7.5cm]{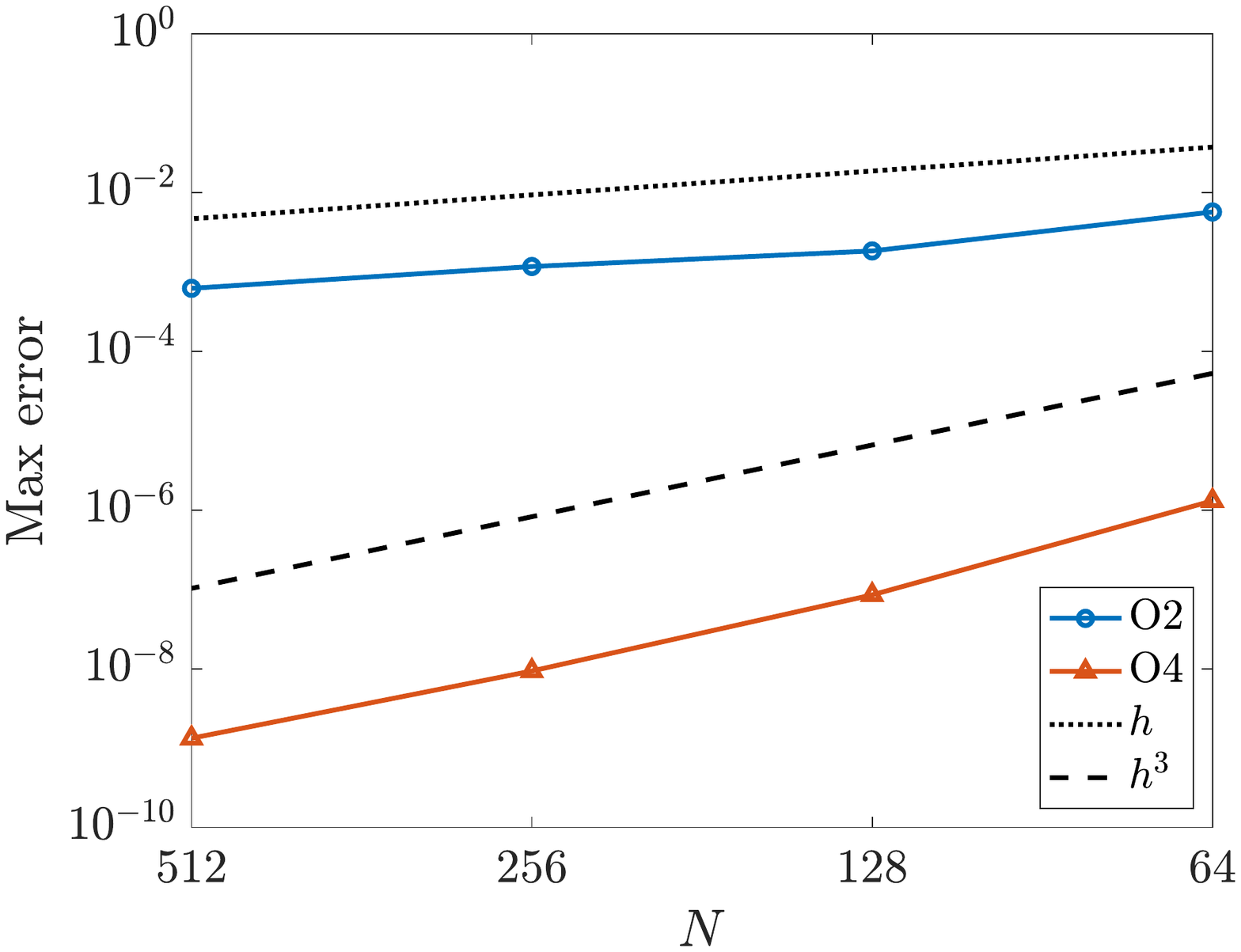}
        \caption{Convergence}
        \label{fig:ellipse_neummann_conv}
    \end{subfigure}
    \quad
    \begin{subfigure}{0.45\textwidth}
        \centering
        \includegraphics[width=\textwidth, trim = 2cm 7cm 2cm 7.5cm]{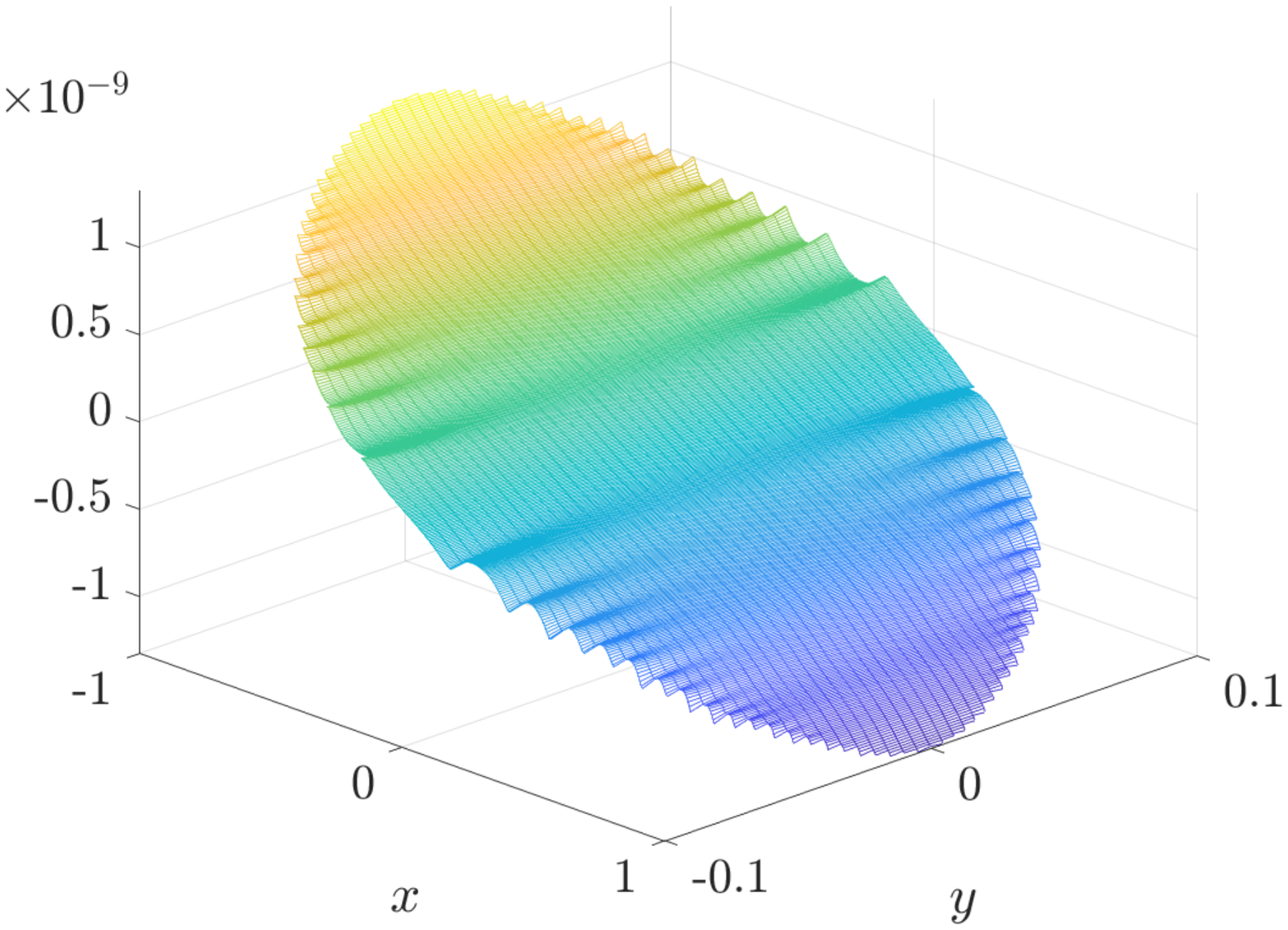}
        \caption{O4 errors $(N=512)$}
        \label{fig:ellipse_neumann_err}
    \end{subfigure}
    \caption{Convergence and errors}\label{fig:ellipse_neumann_conv_err}
\end{figure}

In this section, we look at the same PDE with $\sigma=0$ and the same ellipse with $\alpha=10$ in the same auxiliary domain $[-1.2,1.2]\times[-1.2,1.2]$, with Robin boundary condition
\begin{align}
\frac{\partial u}{\partial n} + u = g
\end{align} 
prescribed at the boundary. Collocation method for boundary method involving normal derivatives will at least lose one order of accuracy, which is expected and can be observed in the convergence plot in Figure~\ref{fig:ellipse_neumann_conv_err}. The errors are smoother near the boundary, compared to the case with Dirichlet boundary conditions. It still gives higher mode oscillations, which we believe the filtering technique in \cite{docampo2020enhancing} can be used to recover better accuracy. Trace finite element \cite{reusken2015analysis} or similar methods targeting surface PDEs can be employed to handle the normal derivatives at the boundary, as the finite element space is already constructed for the boundary element or the cut cell. We will pursue this direction for generalized boundary conditions such as those prescribed by PDEs in future work.

\subsection{A multi-connected domain}
The multi-connected domain used in this example is bounded by the 0 level set of function $\psi(x,y)=\psi_1(x,y)\psi_2(x,y)\psi_3(x,y)$ where the three components are defined as follows
\begin{subequations}
\begin{align}
    \psi_1(x,y) &=  x^2+y^2-1\\
    \psi_2(x,y) &= 4x^2+4(y-0.5)^2-1\\
    \psi_3(x,y) &= 16(x+0.3)^2+16(y+0.4)^2-1
\end{align}
\end{subequations}
The auxiliary domain for this example is $[-1.15,1.15]\times[-1.15,1.15]$.

\begin{figure}[H]
\centering
\begin{tikzpicture}[scale=2]
\draw[thick,fill=gray!40,fill opacity=.5] (0,0) circle (1cm);
\draw[thick,fill=white] (0,0.5) circle (0.5cm);
\draw[thick,fill=white] (-0.3,-0.4) circle (0.25cm);
\draw[->,>=latex] (-1.2,0)--(1.4,0) node[above] {$x$};
\draw[->,>=latex] (0,-1.2)--(0,1.4) node[left] {$y$};
\end{tikzpicture}
\caption{A multi-connected domain}
\label{fig:annulus}
\end{figure}
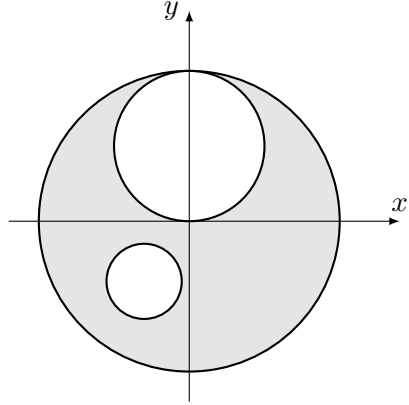

Though there are three separate boundaries in the multi-connected domain, the discrete grid boundary $\gamma$ is still constructed from the stencils of the interior and exterior points, unlike the global basis approach where piecewise smooth global basis functions need to be defined for each circular component of the boundaries, which might become cumbersome when one has more holes, for example, in the application of rigid particles in Stokes flow. 
Second and fourth accurate convergence rates are observed and errors behave similarly as in single-connected domains. The condition numbers are benign and do not grow faster than $1/h$, expect for the first two coarse meshes.



\begin{figure}[H]
    \centering
    \begin{subfigure}{0.45\textwidth}
        \centering
        \includegraphics[width=\textwidth, trim = 2cm 7cm 2cm 7.5cm]{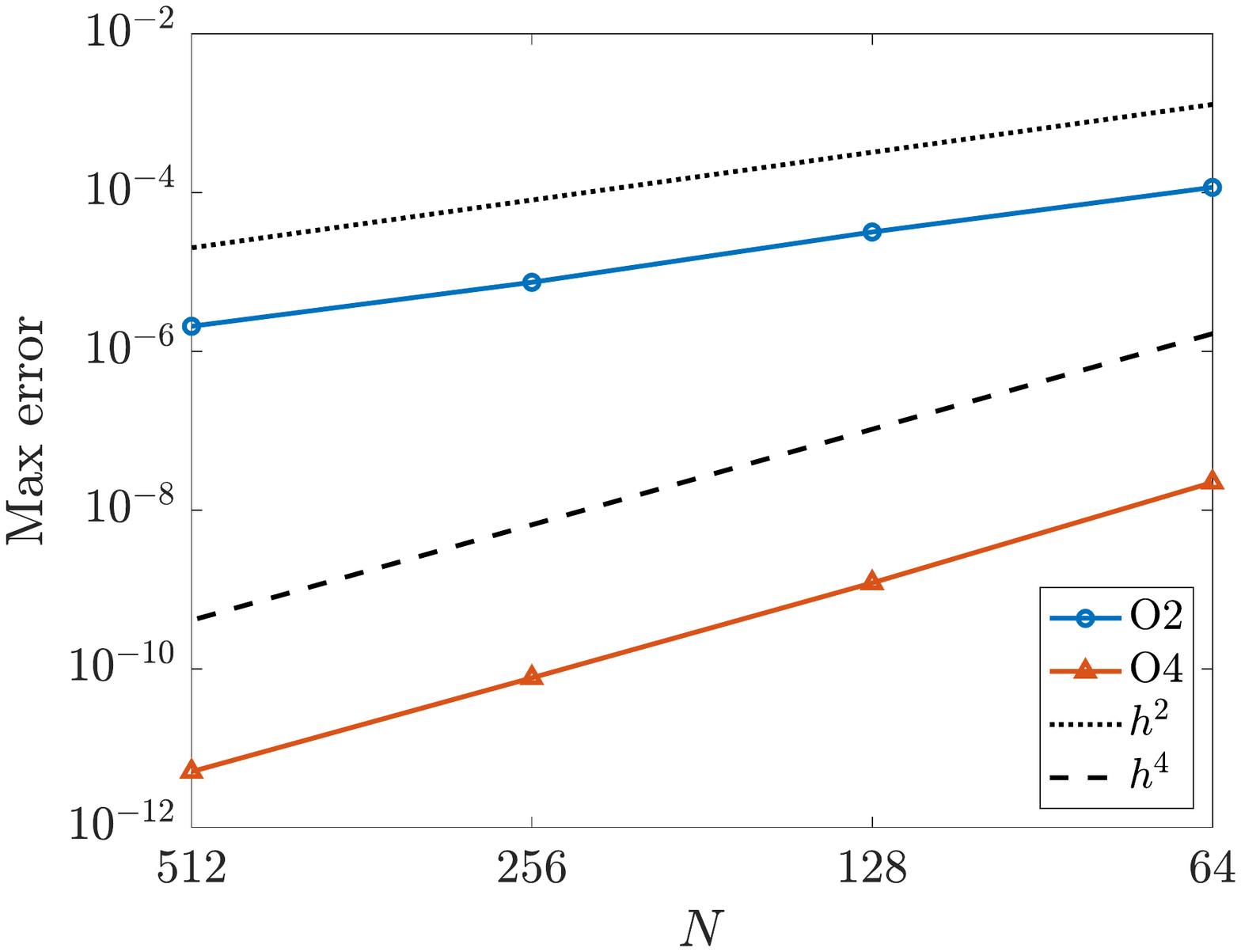}
        \caption{Convergence}
        \label{fig:annulus_conv}
    \end{subfigure}
    \quad
    \begin{subfigure}{0.45\textwidth}
        \centering
        \includegraphics[width=\textwidth, trim = 2cm 7cm 2cm 7.5cm]{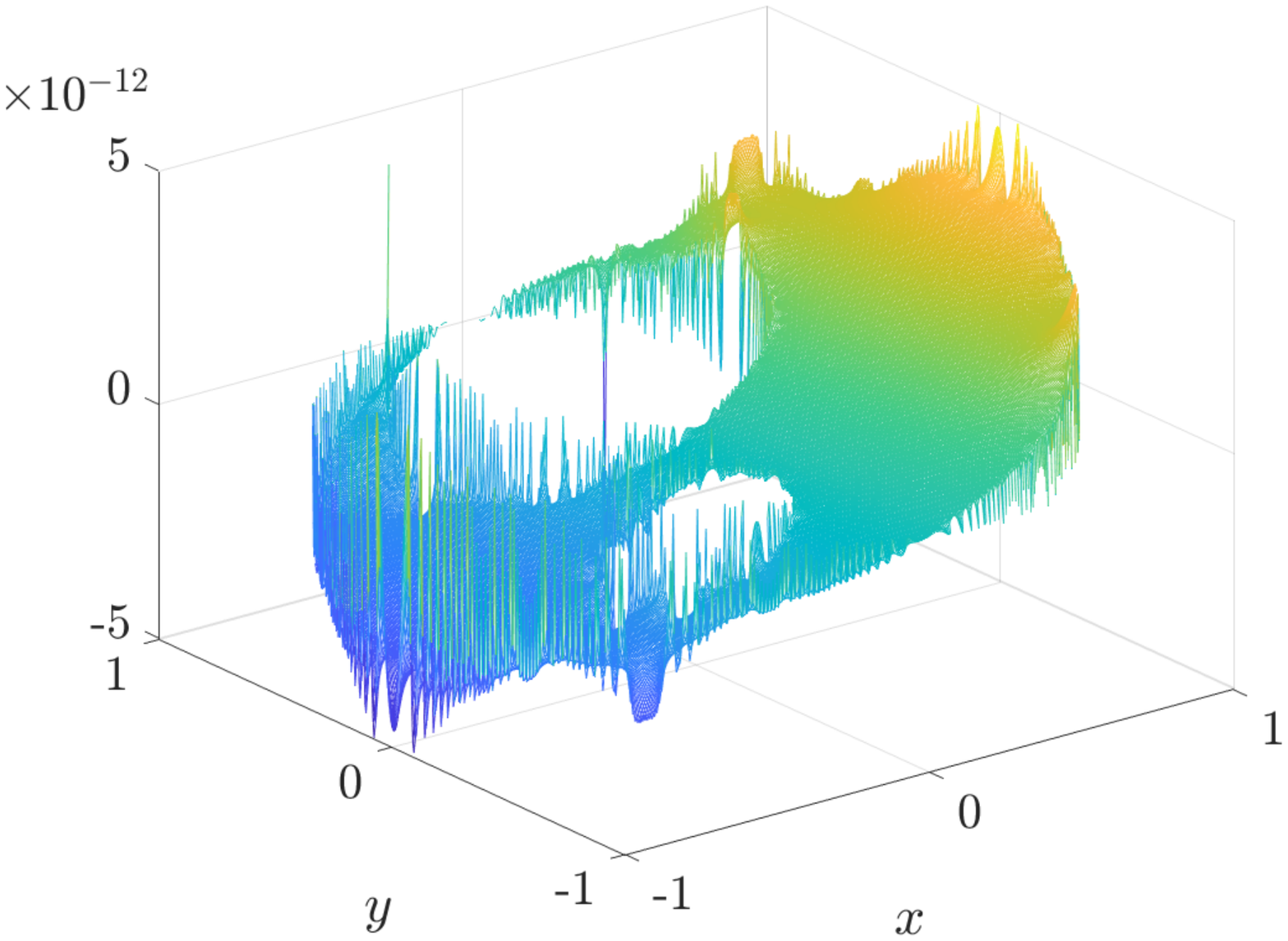}
        \caption{O4 errors $(N=512)$}
        \label{fig:annulus_err}
    \end{subfigure}
    \caption{Convergence and errors}\label{fig:annulus_conv_err}
\end{figure}

\begin{figure}[H]
    \centering
    \includegraphics[width=0.5\textwidth, trim = 2cm 7cm 2cm 7.5cm]{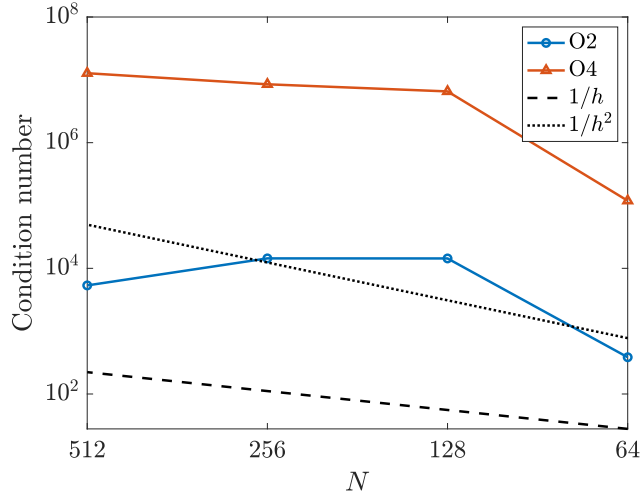}
    \caption{Condition numbers ($\sigma=0$)}
    \label{fig:annulus_cond}
\end{figure}


\subsection{A triangular domain}
In this part, we validate the developed numerical schemes in a triangular domain to show the case of geometry with corners. The triangle has vertices at $(0.5,0.9),(0.9,-0.2),(-0.9,-0.9)$, and is given by the 0 level set of function 
\begin{align}
\psi(x,y) &= -\min(\psi_1(x,y),\psi_2(x,y),1-\psi_1(x,y)-\psi_2(x,y))
\end{align}
where 
\begin{subequations}
\begin{align}
            \psi_1(x,y) &= ((y_3-y_1)(x-x_1)-(x_3-x_1)(y-y_1))/\Delta\\
            \psi_2(x,y) &= ((x_2-x_1)(y-y_1)-(y_2-y_1)(x-x_1))/\Delta
\end{align}
\end{subequations}
with $\Delta$ given by
\begin{align}
\Delta =\left|
\begin{matrix}
1 & x_1 & y_1\\
1 & x_2 & y_2\\
1 & x_3 & y_3\\
\end{matrix}
\right|
\end{align}

The triangle used in this example is given in the following Figure~\ref{fig:triangle} and the auxiliary domain is $[-1.1,1.1]\times[-1.1,1.1]$.
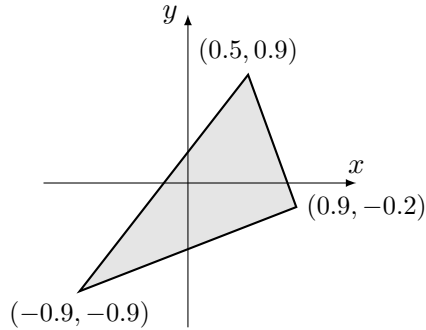
\begin{figure}[H]
\centering
\begin{tikzpicture}[scale=1.6]
\draw[thick,fill=gray!40,fill opacity=.5] (0.5,0.9)
                                        --(0.9,-0.2) 
                                        --(-0.9,-0.9)
                                        -- cycle;
\draw[->,>=latex] (-1.2,0)--(1.4,0) node[above] {$x$};
\draw[->,>=latex] (0,-1.2)--(0,1.4) node[left] {$y$};
\draw(0.5,0.9) node[anchor=south]{\small $(0.5,0.9)$};
\draw(0.9,-0.2)  node[anchor=west]{\small $(0.9,-0.2)$};
\draw(-0.9,-0.9) node[anchor=north]{\small $(-0.9,-0.9)$};
\end{tikzpicture}
\caption{A triangular domain}\label{fig:triangle}
\end{figure}

\begin{figure}[H]
    \centering
    \begin{subfigure}{0.45\textwidth}
        \centering
        \includegraphics[width=\textwidth, trim = 2cm 7cm 2cm 7.5cm]{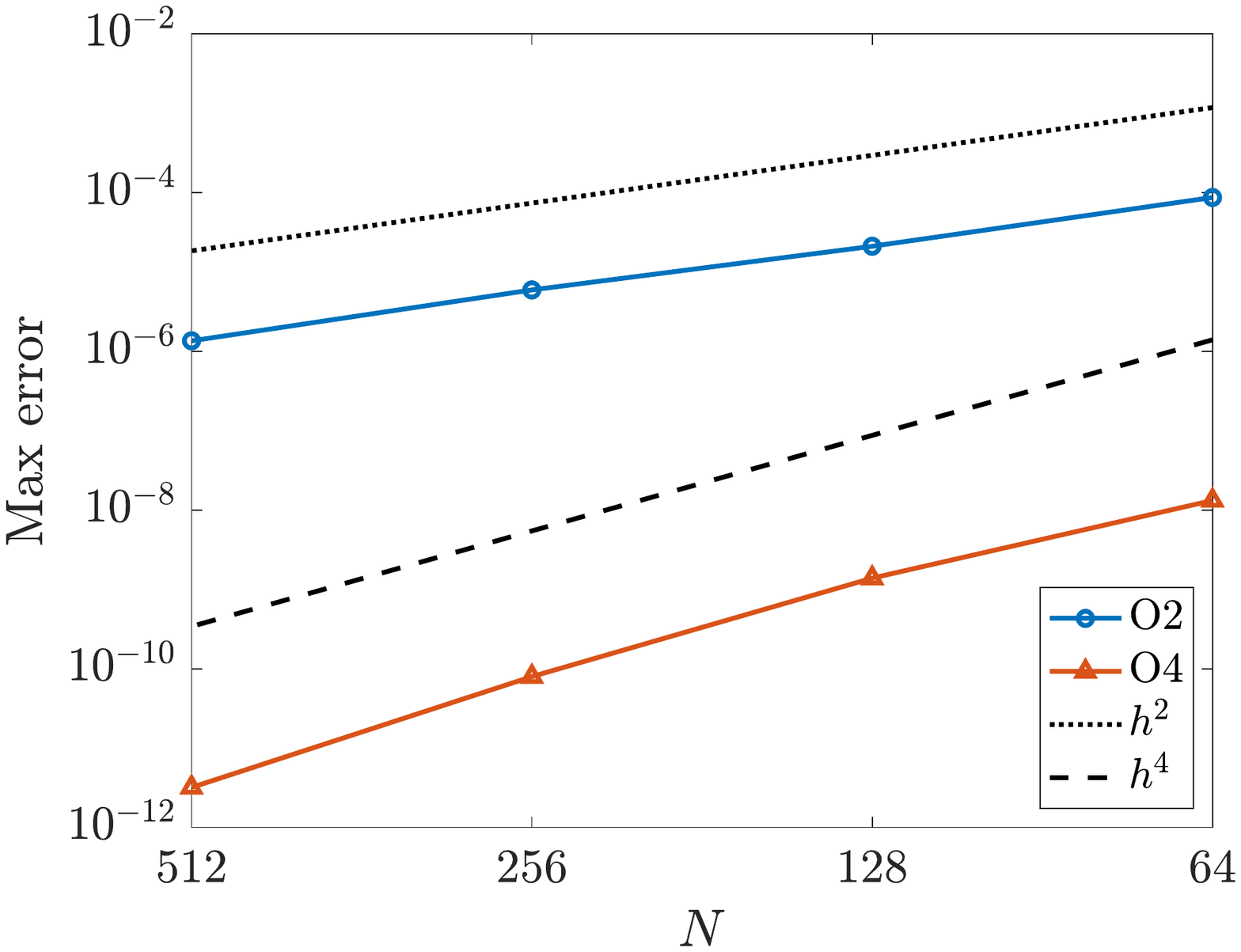}
        \caption{Convergence}
        \label{fig:triangle_conv}
    \end{subfigure}
    \quad
    \begin{subfigure}{0.45\textwidth}
        \centering
        \includegraphics[width=\textwidth, trim = 2cm 7cm 2cm 7.5cm]{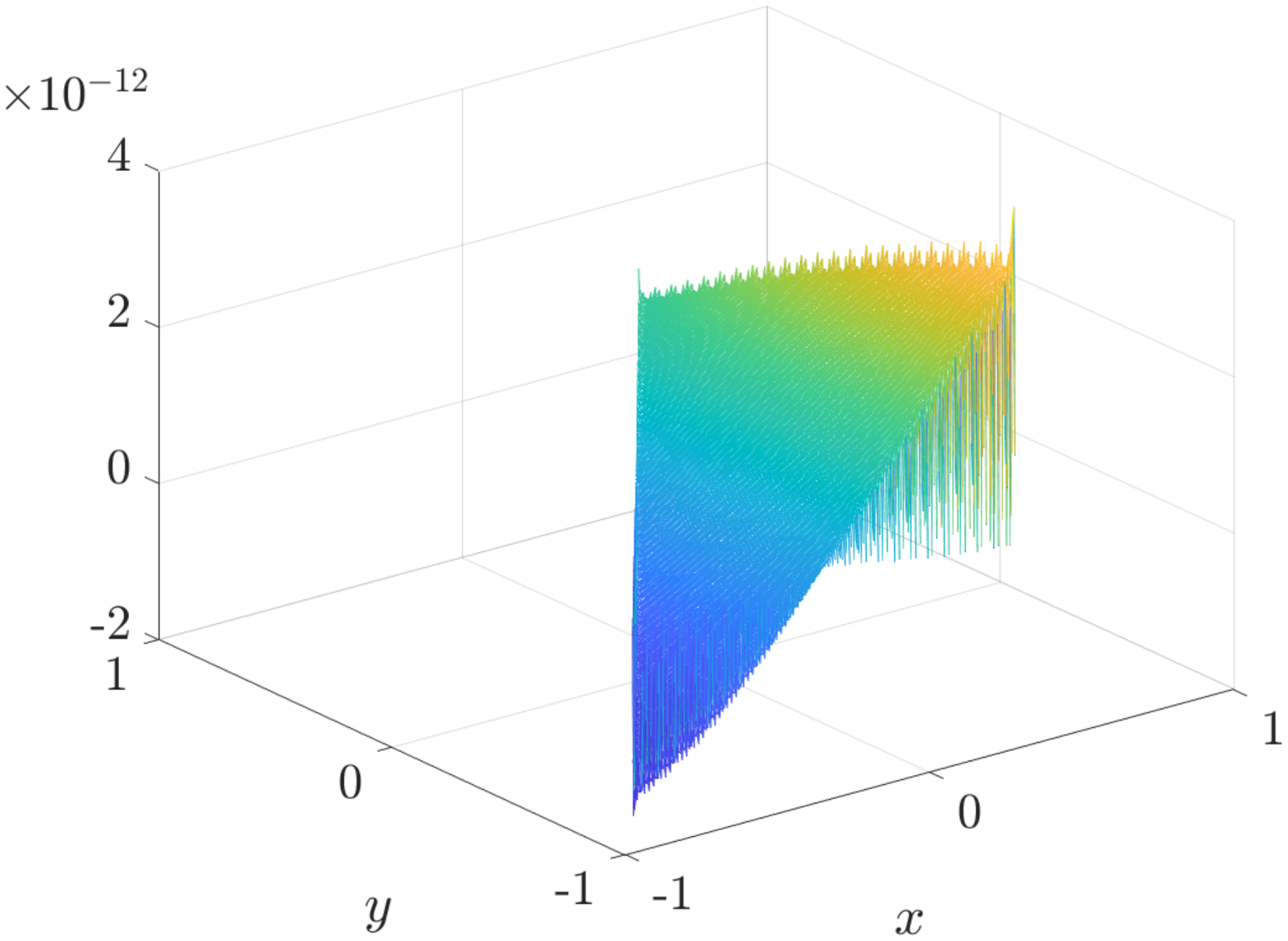}
        \caption{O4 errors $(N=512)$}
        \label{fig:triangle_err}
    \end{subfigure}
    \caption{Convergence and errors}\label{fig:triangle_conv_err}
\end{figure}
Other shapes with corners such as rectangles, diamonds have also been tested, and we only include the example of a triangle to illustrate the robustness of the developed methods.

\begin{figure}[H]
    \centering
    \includegraphics[width=0.5\textwidth, trim = 2cm 7cm 2cm 7.5cm]{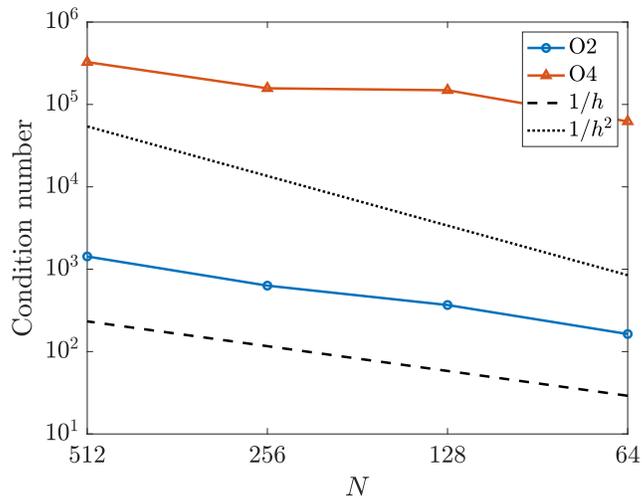}
    \caption{Condition numbers ($\sigma=0$)}
    \label{fig:triangle_cond}
\end{figure}

\section{Conclusions}\label{sec:conclusion}

We developed a geometry-flexible finite difference method within the difference potentials framework using local basis functions. It inherits many features from the classical difference potentials method, such as no need of knowledge of fundamental solutions, the capability of employing fast solvers on simple auxiliary domains, etc. It also brings new features for treating elliptic problems in arbitrary geometry by using local basis functions only at near-boundary grid points. The geometry can be either implicitly or explicitly defined, as we currently only need information of the boundary via the intersection with grid lines. The shapes of geometry can be smooth, piecewise-smooth, multi-connected, disconnected, with corners, etc. We proved the convergence of difference potentials in max norm and validated the developed schemes numerically with different geometries, PDEs and boundary conditions, which also suggests robustness of the developed schemes.

There are of course many places in the current schemes that can be improved. An imminent next step is to consider boundary value problems that involves normal derivatives, where loss of accuracy occurs due to the employment of Lagrange collocation method. This drawback can be overcome by adopting either post-processing techniques or using other boundary treatment rather than collocation, which shall be investigated in future work. The condition numbers of the boundary equations where we solve for the density $u_\gamma$ are always benign and follow the growth $1/h$ for O2. However, in O4 cases they do not share patterns for different geometries and PDEs and can grow faster than $1/h^2$ for certain applications. We will investigate this issue as they are important when iterative solvers are needed.
When it comes to higher dimensions, the developed schemes, in theory, can be extended to three dimensions easily, but the computational cost would be significantly higher. The bottleneck of the computational cost lies within the computations of difference potentials from the multiple unit densities \eqref{eqn:unit_density}. Currently this was partially addressed by using the FFT, which is sufficient for 2D problems. In future applications, we plan to explore GPU acceleration techniques for speed up of computations of difference potentials in 3D, as the difference potentials associated with different unit densities are independent, which means the computations are highly parallelizable. It can also be argued that the computational cost of difference potentials can be hid in precomputation for fixed geometry and time dependent problems. It is one of our future research directions to design an efficient 3D version of the the algorithms in this work. Other possible future research directions include interface problems such as in composite materials, bulk-surface problems, nonlinear problems, etc.

\section*{Acknowledgement}
The author would like to thank COMSOL for funding. The author is grateful to Jeffrey Banks, Elias Jarlebrin, Jennifer Ryan, Anna-Karin Tornberg and the group, Siyang Wang and Sara Zahedi for many helpful discussions.




\bibliographystyle{elsart-num}
\bibliography{ref.bib}
 
\end{document}